\documentclass[11pt]{amsart}

\usepackage[margin=0.9in]{geometry}
\usepackage[foot]{amsaddr}

\usepackage[T1]{fontenc}
\usepackage{lmodern}

\usepackage{setspace}
\usepackage{amssymb}
\usepackage{amsmath}
\usepackage{latexsym}
\usepackage[scr=boondoxo,scrscaled=1.05]{mathalfa}
\usepackage{tensor}
\usepackage{upgreek}
\usepackage{graphicx}
\usepackage{lscape}
\usepackage{booktabs}

\theoremstyle{remark}
\newtheorem{remark}{Remark}

\begin{document}

\title[A Fast-Convergence Resolution of the Stochastic Eigenproblem]{A Fast-Convergence Resolution of the Stochastic Eigenproblem Using Halley's Method and the Spectral-Chaos Approach}

\author{Hugo Esquivel}
\address{Lead Structural Engineer, Estructuras del Norte, Barranquilla, Colombia}
\email{hesquivel-@estrunorte.com.co}

\author{Kabir Oluwatobi Idowu}
\address{Department of Mathematics, Purdue University, West Lafayette, Indiana, 47907}
\email{kidowu@purdue.edu}

\author{Guang Lin}
\address{Department of Mathematics, School of Mechanical Engineering, Purdue University, West Lafayette, Indiana, 47907}
\email{guanglin@purdue.edu}

\maketitle

\begin{abstract}
\noindent Solving stochastic eigenvalue problems has long been essential for informed decision-making, advancing scientific knowledge, and ensuring the reliability of engineering designs and applications.
This paper underscores the need to continue enhancing existing numerical methods for solving the stochastic eigenproblem in order to improve convergence rates, computational efficiency, and robustness.
Specifically, we propose a novel spectral-chaos method for solving the stochastic (linear) eigenvalue problem, employing Halley's method as the root-finding algorithm to leverage its cubic convergence properties.
Our method achieves maximal convergence in solving stochastic eigenvalue problems since its rate cannot be further improved using a higher-order Householder method due to the quadratic nature of the resulting system of equations.
Additionally, due to the complexity of the resulting system of equations, a tensorial approach was developed to tackle the challenges associated with the dimensional multiplicity of the stochastic eigenvalue problem, without which the solution would have been intractable.
The method is derived rigorously, with a detailed error analysis that highlights the benefit of using our approach when the eigenvector components are nearly known, the computational cost of the method is also rigorously presented, and an illustrative example is provided to demonstrate the implementation of the method.
Subsequently, a case study is demoed to analyze the results and validate the advantages of using Halley's method over Newton's method and Monte Carlo simulations.
\end{abstract}

\section{Introduction}

The complexity of stochastic processes and their associated eigenvalue problems presents ongoing challenges in computational mathematics and applied sciences \cite{Frank2024SpectralEigenvalues,Lotz2020WilkinsonsEigenproblems,esquivel2021efficient}.
Solving these eigenproblems is critical for applications in quantum physics, mathematics, engineering, and financial modeling, where stochastic models effectively represent inherent uncertainties \cite{Eibert2012SolvingEigenvalue,khan2024optimization}.
However, conventional methods often face limitations in effectiveness and accuracy, particularly when applied to complex high-dimensional systems \cite{lucas2018using, agarwala2024high}.
These challenges underscore the need for robust and efficient methods to address the stochastic eigenvalue problem in a manner that balances computational complexity with solution accuracy.

In this paper, we study the eigenspaces associated with stochastically governed linear dynamical models.
These models act as a substitute for the observed random spatial variations in the mechanical representation of the system \cite{li2009stochastic}.
Such representations are crucial for capturing the probabilistic behavior of dynamical systems under uncertainty, particularly in engineering and physical sciences, where precise predictions of modal quantities are essential.
The resulting stochastic eigenvalue problem is critical for understanding dynamical systems as a whole \cite{adhikari2007random,san2000stochastic}, as statistical measures of modal quantities have a substantial impact on the expected performance of these systems.
For example, in structural engineering, stochastic eigenvalue analyses guide the design of structures to ensure the stability and reliability of the system under uncertain loading conditions \cite{aldosary2018structural,antonio2017reliability,lee2009stochastic}.
In this case, stochastic reduced order models are often more successful because they effectively represent the observed variability of the system \cite{mou2023efficient,frohlich2022uncertainty,simiu2009chaotic}.
Historically, these problems have been addressed using low-order Taylor expansions with perturbation methods or by analyzing statistical properties using low-order statistics \cite{pryse2021neumann,charalampopoulos2022uncertainty}.
However, these approaches often fail to capture higher-order interactions, especially in systems with significant nonlinearity or higher-dimensional randomness \cite{gu2022change,medina2021hyperharmonic}.

In contemporary science and engineering, random matrices are increasingly used as random perturbations of finite-dimensional operators to model uncertainty \cite{ding2025advanced,matthies2008stochastic,schenk2005uncertainty}.
Unlike matrices derived from finite-dimensional representations of partial differential operators, these random matrices often allow closed-form expressions of the statistical moments and probability density functions of their eigenvalues and eigenvectors in specific scenarios \cite{couillet2022random,bai2010spectral}.
However, the practical applicability of such methods is limited, particularly for large-scale systems or those characterized by complex parameter uncertainties.
Conversely, the matrices discussed in this paper arise from finite-dimensional approximations of continuous systems, where randomness stems from uncertainties in the system's parameters.
For these systems, closed-form solutions to the random eigenvalue problem are typically unavailable, necessitating alternative strategies such as statistical sampling \cite{adhikari2007random,vidyasagar1998statistical}, perturbation techniques \cite{huang2018homotopy,bold2006heterogeneous}, or polynomial chaos representations combined with Galerkin projections \cite{du2025uncertainty,wang2021interval,arnst2014reduced}.
Although statistical sampling-based methods are conceptually straightforward, they are often computationally demanding and exhibit poor scalability with increasing system complexity \cite{ruthotto2021introduction,kougioumtzoglou2020sparse}.

Halley's method, known for its cubic convergence rate in root-finding problems, offers a promising approach to achieving reliable solutions \cite{Basto2017ConvergenceMethods,Eustaquio2018AClass,Cordero2020SomeProblems}.
This superlinear convergence property makes it particularly suitable for solving nonlinear systems that exhibit rapid changes in solution behavior \cite{algelany2022chaotic,iyengar2025polynomial}.
However, its application to stochastic eigenproblems remains underexplored.
For example, \cite{ghanem2007efficient,leader2022numerical} used polynomial chaos expansions to transform the stochastic problem into a deterministic one, thus allowing the computation of eigenvalues as functions of random variables.
Although effective, the quadratic rate of convergence inherent in these methods imposes limitations on computational efficiency, particularly for problems with high degrees of randomness or complex eigenvalue spectra.
Furthermore, reliance on deterministic transformations can lead to loss of precision in capturing intricate stochastic variations, especially in cases involving non-Gaussian randomness or multimodal distributions \cite{lu2015limitations,ghanem2017polynomial}.

This paper addresses these challenges by introducing an innovative application of Halley's method within the spectral-chaos framework for solving stochastic eigenvalue problems.
Halley's method, also known as the second-order Householder method, provides a cubic convergence rate---a significant improvement over the quadratic rate offered by Newton-based methods traditionally employed in polynomial chaos expansions \cite{lloyd2022fast,obsieger2013numerical,lloyd2023precision}.
By reformulating the stochastic eigenvalue problem to enable the application of Halley's method in a multidimensional random-spatial setting, the proposed approach achieves a polynomial-order enhancement in convergence rate and the ability to accelerate when the eigenvector components are nearly known.
This innovation preserves the strengths of polynomial chaos while substantially reducing computational effort and increasing precision.
Additionally, the introduction of a tensorial form of the stochastic operator facilitates the efficient manipulation of the random-spatial space, making this approach particularly effective for addressing the dimensional multiplicity of the stochastic eigenvalue problem.

Beyond its theoretical contributions, the proposed method offers practical benefits for computational efficiency and accuracy.
For example, the improved convergence rate significantly reduces the number of iterations required to achieve high precision, directly addressing the computational bottlenecks of existing methods.
Moreover, the ability to capture higher-order stochastic interactions makes our improved spectral-chaos approach suitable for a broader range of problems, including those involving nonstandard distributions or parameter uncertainties.

The following sections delve into the details of Halley's method for solving stochastic (linear) eigenvalue problems and its implications for linear dynamical systems.
Here is a brief outline of the contents of each section.
Section \ref{sec:prosta} outlines the problem statement.
Section \ref{sec:spesol} introduces the spectral-chaos solution framework.
Section \ref{sec:promom} explores the probability moments of the eigenpairs.
Section \ref{sec:numsol} presents the numerical scheme used to solve the resulting system of quadratic equations with Halley's method.
Section \ref{sec:errana} gives a detailed error analysis for the proposed method.
Section \ref{sec:comcos} provides the computational cost of the proposed method.
Section \ref{sec:illexa} demonstrates the method through an illustrative example.
Section \ref{sec:casstu} demos a case study to validate the robustness and efficiency of the proposed approach.
Finally, Section \ref{sec:conclu} offers concluding remarks.

\section{Problem statement}\label{sec:prosta}

Before stating the problem, it is necessary to establish a series of definitions.

Let $(\Theta,\boldsymbol{\Theta},\nu)$ be a \emph{probability space}, where $\Theta$ is the sample space (the set of all possible outcomes), $\boldsymbol{\Theta}\subset 2^\Theta$ is the $\sigma$-algebra on $\Theta$ (the collection of events), and $\nu:\boldsymbol{\Theta}\to[0,1]$ is the probability measure on $\boldsymbol{\Theta}$.

Let $(\mathfrak{Z},\boldsymbol{\mathfrak{Z}},\mu)$ be the \emph{random space} related to $(\Theta,\boldsymbol{\Theta},\nu)$ via the measurable function $\zeta:(\Theta,\boldsymbol{\Theta})\to(\mathbb{R}^d,\mathcal{B}_{\mathbb{R}^d})$ defined by $\zeta=\zeta(\omega)$, where $\mathfrak{Z}=\zeta(\Theta)\subset\mathbb{R}^d$ is the random domain (that is, the $\mathbb{R}^d$ representation of $\Theta$), $\boldsymbol{\mathfrak{Z}}=\mathcal{B}_{\mathbb{R}^d}\cap \mathfrak{Z}$ is the $\sigma$-algebra on $\mathfrak{Z}$, $\mu:\boldsymbol{\mathfrak{Z}}\to[0,1]$ is the probability measure on $\boldsymbol{\mathfrak{Z}}$ defined by the pushforward of $\nu$ by $\zeta$ (that is, $\mu=\zeta_*(\nu)$), $\mathcal{B}_{\mathbb{R}^d}$ is the Borel $\sigma$-algebra on $\mathbb{R}^d$, and $d$ represents the dimensionality of the random space.

Let $(\mathfrak{E},\boldsymbol{\eta})$ be the \emph{spatial space}, where $\mathfrak{E}=\mathbb{R}^R$ is the spatial domain, $R$ denotes the dimensionality of the space, and $\boldsymbol{\eta}$ is a Riemannian metric tensor (as specified below).
With this metric tensor, $\mathfrak{E}$ is a Riemannian manifold.

Let $\mathscr{Z}=L^2(\mathfrak{Z},\boldsymbol{\mathfrak{Z}},\mu; \mathbb{R})$ be a Lebesgue square-integrable space equipped with its standard inner product
\begin{equation*}
\langle\cdot\,,\cdot\rangle:\mathscr{Z}^2\to \mathbb{R}
\quad:\Leftrightarrow\quad
(f,g)\mapsto\langle f,g\rangle=\int fg\,\mathrm{d}\mu.
\end{equation*}
This space is referred to as the \emph{random function space}, and it is a Hilbert space.

Finally, let $\mathrm{L}(n,\mathbb{R})$ be the space of linear maps from $\mathbb{R}^n$ to $\mathbb{R}^n$, $\mathrm{Sym}(n,\mathbb{R})\subset\mathrm{L}(n,\mathbb{R})$ the space of symmetric linear maps from $\mathbb{R}^n$ to $\mathbb{R}^n$, and $\mathrm{GL}(n,\mathbb{R})\subset\mathrm{L}(n,\mathbb{R})$ the space of linear isomorphisms from $\mathbb{R}^n$ to $\mathbb{R}^n$.

The problem is to find the $n$-th random eigenpair $(\lambda,\boldsymbol{\upphi}):\mathfrak{Z}\to\mathbb{R}^{R+1}$ with $n\in\{1,2,\ldots,R\}$, such that ($\mu$-a.e.)
\begin{subequations}\label{eq:prosta100}
\begin{align}
\mathbf{K}\boldsymbol{\upphi}=\lambda\boldsymbol{\upphi}&\quad\text{on $\mathfrak{Z}$}\label{eq:prosta100a}\\
\boldsymbol{\upphi}\boldsymbol{\cdot}\boldsymbol{\upphi}=1&\quad\text{on $\mathfrak{Z}$},\label{eq:prosta100b}
\end{align}
\end{subequations}
where $\mathbf{K}:\mathfrak{Z}\to\mathrm{Sym}(R,\mathbb{R})$ is a random, real symmetric matrix of size $R$, $\lambda:\mathfrak{Z}\to\mathbb{R}$ is the $n$-th random eigenvalue, $\boldsymbol{\upphi}:\mathfrak{Z}\to\mathfrak{E}$ is the $n$-th random eigenvector, and $\,\boldsymbol{\cdot}\,:\mathfrak{E}^2\to\mathbb{R}$ is the dot product in linear algebra.
In \cite{ghanem2007efficient} it was shown that if all entries of $\mathbf{K}$ are in $\mathscr{Z}$, so is $\lambda$, and if Eq.~\eqref{eq:prosta100b} is considered to hold $\mu$-a.e., it is because all entries of $\boldsymbol{\upphi}$ are in $\mathscr{Z}$.

More generally, Eq.~\eqref{eq:prosta100} can be expressed in tensor form as ($\mu$-a.e.)
\begin{subequations}\label{eq:prosta200}
\begin{align}
\boldsymbol{K}(\zeta)[\,\cdot\,,\boldsymbol{\phi}(\zeta)]=\lambda(\zeta)\,\boldsymbol{\phi}(\zeta)\quad\forall\zeta\in\mathfrak{Z} & \quad\Leftrightarrow\quad K\indices{^u_v}\phi^v=\lambda\phi^u\label{eq:prosta200a}\\
\boldsymbol{\eta}[\boldsymbol{\phi}(\zeta),\boldsymbol{\phi}(\zeta)]=1\quad\forall\zeta\in\mathfrak{Z} & \quad\Leftrightarrow\quad\eta\indices{_{uv}}\phi^u\phi^v=1,\label{eq:prosta200b}
\end{align}
\end{subequations}
where $\boldsymbol{K}(\zeta)=K\indices{^u_v}(\zeta)\,\mathbf{e}_u\otimes\mathbf{e}^v:\mathfrak{E}^*\times\mathfrak{E}\to\mathbb{R}$ is a spatial $(1,1)$-tensor with the property that its flat, $\boldsymbol{K}^\flat(\zeta)$, is symmetric, $\boldsymbol{\phi}(\zeta)=\phi^u(\zeta)\,\mathbf{e}_u$ is a spatial vector, $\mathbf{e}_u\in\mathfrak{E}$ is the standard unit vector in the $u$-th dimension, $\mathbf{e}^v:\mathfrak{E}\to\mathbb{R}$ is the standard covector in the $v$-th dimension, $\boldsymbol{\eta}=\eta_{uv}\,\mathbf{e}^u\otimes\mathbf{e}^v:\mathfrak{E}^2\to\mathbb{R}$ is a Riemannian metric tensor endowed on $\mathfrak{E}$ (typically, the Euclidean metric tensor), $\mathfrak{E}^*$ is the dual space of $\mathfrak{E}$ with $\mathbf{e}^v[\mathbf{e}_u]=\delta^v_u$ (the Kronecker delta), $K\indices{^u_v},\phi^u,\lambda:\mathfrak{Z}\to\mathbb{R}$ are random functions in $\mathscr{Z}$, and $u,v\in\{1,2,\ldots,R\}$.
As is customary in multi-linear algebra, a summation is always implied over repeated indices---this is known as the Einstein summation convention.

\begin{remark}
A more compelling way to define the tensor field $\boldsymbol{K}$ described above is by defining it as a $\mathscr{Z}$ section of $\mathrm{T}\indices{^1_1}(\mathbf{Z}):=\mathbf{Z}\otimes\mathbf{Z}^*$, where $\mathbf{Z}:=\mathfrak{Z}\times\mathbb{R}^R\xrightarrow{\pi}\mathfrak{Z}$ is the trivial vector bundle of rank $R$ over $\mathfrak{Z}$.
This can be accomplished through a two-step process as follows.
First, let:
\begin{equation*}
\boldsymbol{K}:\mathfrak{Z}\to\mathrm{T\indices{^1_1}}(\mathbf{Z})
\quad:\Leftrightarrow\quad
\zeta\mapsto\boldsymbol{K}(\zeta)=K\indices{^u_v}(\zeta)\,\mathbf{e}_u\otimes\mathbf{e}^v\in\mathrm{T}\indices{^1_1}(\mathbf{Z}_\zeta)
\quad\text{with}\quad
K\indices{^u_v}\in\mathscr{Z},
\end{equation*}
where $\mathbf{Z}_\zeta=\pi^{-1}(\{\zeta\})=\{\zeta\}\times\mathbb{R}^R$ is the fiber over $\zeta\in\mathfrak{Z}$, which can be identified with $\mathfrak{E}=\mathbb{R}^R$.
Then, define:
\begin{equation*}
\boldsymbol{K}(\zeta)=K\indices{^u_v}(\zeta)\,\mathbf{e}_u\otimes\mathbf{e}^v:(\mathbf{Z}_\zeta)^*\times\mathbf{Z}_\zeta\to\mathbb{R}
\quad:\Leftrightarrow\quad
(\boldsymbol{\alpha},\mathbf{w})\mapsto\boldsymbol{K}(\zeta)[\boldsymbol{\alpha},\mathbf{w}]=\alpha_u w^v\, K\indices{^u_v}(\zeta),
\end{equation*}
where $\boldsymbol{\alpha}=\alpha_u\mathbf{e}^u$, and $\mathbf{w}=w^v\mathbf{e}_v$.
Put differently, $\boldsymbol{K}$ is a $\mathscr{Z}$-tensor field of order $(1,1)$ on $\mathfrak{Z}$, assigning to each point $\zeta\in\mathfrak{Z}$ an element $\boldsymbol{K}(\zeta)$ of $\mathrm{T}\indices{^1_1}(\mathbf{Z}_\zeta)\cong\mathfrak{E}\otimes\mathfrak{E}^*$.
However, for the sake of simplicity, we opted for the more straightforward, albeit more decorated, definition presented earlier.
Likewise, $\boldsymbol{\phi}$ is a $\mathscr{Z}$ section of $\mathbf{Z}$ (or a $\mathscr{Z}$-vector field on $\mathfrak{Z}$), and $\boldsymbol{\eta}=\boldsymbol{\eta}[\,\cdot\,,\cdot\,]$ can be reinterpreted as a $\mathscr{C}$ section of $\mathrm{T}\indices{^0_2}(\mathbf{Z})=\mathbf{Z}^*\otimes\mathbf{Z}^*$, with $\mathscr{C}\subset\mathscr{Z}$ symbolizing the space of real-valued constant functions on $\mathfrak{Z}$; hence, admitting the form of a tensor field, $\boldsymbol{\eta}=\boldsymbol{\eta}(\cdot)[\,\cdot\,,\cdot\,]$.

\end{remark}

\begin{remark}
Because by definition $\mathbf{K}(\zeta)\in\mathrm{Sym}(R,\mathbb{R})$, we require that the flat of $\boldsymbol{K}(\zeta)$:
\begin{equation*}
\boldsymbol{K}^\flat(\zeta)=\eta_{uw}K\indices{^w_v}(\zeta)\,\mathbf{e}^u\otimes\mathbf{e}^v:\mathfrak{E}^2\to\mathbb{R}
\end{equation*}
is symmetric as well; i.e.,~we require that $\boldsymbol{K}^\flat(\zeta)[\mathbf{x},\mathbf{y}]=\boldsymbol{K}^\flat(\zeta)[\mathbf{y},\mathbf{x}]$ for all $\zeta\in\mathfrak{Z}$ and $\mathbf{x},\mathbf{y}\in\mathfrak{E}$.
\end{remark}

\section{Spectral solution of eigenproblem}\label{sec:spesol}
Let $(\varphi_j:(\mathfrak{Z},\boldsymbol{\mathfrak{Z}})\to(\mathbb{R},\mathcal{B}_\mathbb{R}))_{j=0}^\infty$ be an ordered orthogonal basis in $\mathscr{Z}$ with $\varphi_0\equiv1$.
If $f$ is an element of $\mathscr{Z}=\mathrm{span}(\varphi_j)_{j=0}^\infty$, then
\begin{equation}\label{eq:spesol100}
f=\sum_{j=0}^\infty f^j\varphi_j\equiv f^j\varphi_j,
\end{equation}
where $f^j\in\mathbb{R}$ is the $j$-th component of $f$.
The orthogonality property of the basis implies that
\begin{equation}\label{eq:spesol150}
\langle\varphi_i,\varphi_j\rangle=\langle\varphi_i,\varphi_i\rangle\,\delta^i_j.
\end{equation}

Moreover, define $\mathscr{Z}'$ to be the dual space of $\mathscr{Z}$ spanned by the continuous linear functionals $(\varphi^i:\mathscr{Z}\to\mathbb{R})_{i=0}^\infty$ given by
\begin{equation}\label{eq:spesol200}
\varphi^i[f]:=[\varphi^i,f]=\frac{\langle\varphi_i,f\rangle}{\langle\varphi_i,\varphi_i\rangle}=f^j\frac{\langle\varphi_i,\varphi_j\rangle}{\langle\varphi_i,\varphi_i\rangle}=f^j\delta^i_j=f^i,
\end{equation}
where $[\,\cdot\,,\cdot\,]:\mathscr{Z}'\times\mathscr{Z}\to\mathbb{R}$ is the natural pairing of $\mathscr{Z}$ and $\mathscr{Z}'$ given by the first equality.
The second equality follows from Eq.~\eqref{eq:spesol100}, and the third equality from Eq.~\eqref{eq:spesol150}.
Thus, when $\varphi^i$ acts on $f\in\mathscr{Z}$, it outputs the $i$-th component of $f$.
Another result from the definition is that $\varphi^i[\varphi_j]=\delta^i_j$.

Due to the assumption that both $\lambda$ and $\phi^v$ are in $\mathscr{Z}$, we have
\begin{subequations}\label{eq:spesol300}
\begin{align}
\lambda=\lambda^j\varphi_j:\mathscr{Z}\to\mathbb{R}\quad &:\Leftrightarrow\quad\zeta\mapsto\lambda(\zeta)=\sum_{j=0}^\infty\lambda^j\,\varphi_j(\zeta)\equiv\lambda^j\,\varphi_j(\zeta)\\
\phi^v=\phi^{vj}\varphi_j:\mathscr{Z}\to\mathbb{R}\quad &:\Leftrightarrow\quad\zeta\mapsto\phi^v(\zeta)=\sum_{j=0}^\infty\phi^{vj}\,\varphi_j(\zeta)\equiv\phi^{vj}\,\varphi_j(\zeta).
\end{align}
\end{subequations}

Substituting Eq.~\eqref{eq:spesol300} into Eq.~\eqref{eq:prosta200} gives
\begin{subequations}\label{eq:spesol400}
\begin{gather}
K\indices{^u_v}(\phi^{vj}\varphi_j)=(\lambda^k\varphi_k)(\phi^{uj}\varphi_j)\\
\eta_{uv}(\phi^{uj}\varphi_j)(\phi^{vk}\varphi_k)=\varphi_0\varphi_0,
\end{gather}
\end{subequations}
where $\varphi_0\varphi_0\equiv 1$ because $\varphi_0\equiv 1$ by definition.

Projecting Eq.~\eqref{eq:spesol400} onto $\mathscr{Z}$ yields
\begin{subequations}\label{eq:spesol500}
\begin{gather}
\varphi^i[K\indices{^u_v}\phi^{vj}\varphi_j]=\varphi^i[\lambda^k\phi^{uj}\varphi_j\varphi_k]\\
\varphi^i[\eta_{uv}\phi^{uj}\phi^{vk}\varphi_j\varphi_k]=\varphi^i[\varphi_0\varphi_0],
\end{gather}
\end{subequations}
which after simplification produces a system of infinitely many quadratic equations:
\begin{subequations}\label{eq:spesol600}
\begin{gather}
K\indices{^u_v^i_j}\phi^{vj}=I\indices{^i_{jk}}\lambda^k\phi^{uj}\\
I\indices{^i_{jk}}\eta_{uv}\phi^{uj}\phi^{vk}=I\indices{^i_{00}}
\end{gather}
\end{subequations}
with $K\indices{^u_v^i_j}$ and $I\indices{^i_{jk}}$ given by
\begin{equation}\label{eq:spesol700}
K\indices{^u_v^i_j}=\frac{\langle\varphi_i,K\indices{^u_v}\varphi_j\rangle}{\langle\varphi_i,\varphi_i\rangle}
\quad\text{and}\quad
I\indices{^i_{jk}}=\frac{\langle\varphi_i,\varphi_j\varphi_k\rangle}{\langle\varphi_i,\varphi_i\rangle}.
\end{equation}

\begin{remark}
If the expressions in Eq.~\eqref{eq:spesol700} do not yield real numbers, it means that the problem cannot be solved using the spectral approach.
This is because $\mathscr{Z}$ does not form an algebra.
\end{remark}

\begin{remark}
Eq.~\eqref{eq:spesol700} entails the existence of the random-spatial tensors:
\begin{subequations}\label{eq:spesol800}
\begin{gather}
\boldsymbol{\tilde{K}}=K\indices{^u_v^i_j}\,\mathbf{e}_u\otimes\mathbf{e}^v\otimes\varphi_i\otimes\varphi^j:\mathfrak{E}^*\times\mathfrak{E}\times\mathscr{Z}'\times\mathscr{Z}\to\mathbb{R}\\
\boldsymbol{\tilde{I}}=I\indices{^i_{jk}}\,\varphi_i\otimes\varphi^j\otimes\varphi^k:\mathscr{Z}'\times\mathscr{Z}^2\to\mathbb{R},
\end{gather}
\end{subequations}
where $\boldsymbol{\tilde{I}}$ is symmetric in the indices $j$ and $k$, and the tensor
\begin{equation*}
\boldsymbol{\tilde{K}}^{\flat_1}=K\indices{_{uv}^i_j}\,\mathbf{e}^u\otimes\mathbf{e}^v\otimes\varphi_i\otimes\varphi^j:\mathfrak{E}^2\times\mathscr{Z}'\times\mathscr{Z}\to\mathbb{R}
\end{equation*}
is symmetric in the indices $u$ and $v$ because $\mathbf{K}$ is a symmetric matrix by definition.
Since $\flat_1$ denotes the flat of $\boldsymbol{\tilde{K}}$ with respect to its first slot, we have $K\indices{_{uv}^i_j}=\eta_{uw}K\indices{^w_v^i_j}$.
\end{remark}

\section{Probability moments of eigenpair}\label{sec:promom}

The key probability moments are the mean and variance of the eigenpair $(\lambda,\boldsymbol{\phi})$.
If $\varrho$ denotes the $k$-th coordinate of $(\lambda,\boldsymbol{\phi})$ with $k\in\{1,2,\ldots,R+1\}$, then we have
\begin{equation}\label{eq:promom100}
\varrho(\zeta)=\sum_{j=0}^\infty \varrho^j\,\varphi_j(\zeta)\equiv \varrho^j\,\varphi_j(\zeta),
\end{equation}
where the $j$-th component of $\varrho$ is given by
\begin{equation*}
\varrho^j=\varphi^j[\varrho]=\frac{\langle\varphi_j,\varrho\rangle}{\langle\varphi_j,\varphi_j\rangle}.
\end{equation*}

Moreover, let the standard metric tensor on $\mathscr{Z}$ be taken as
\begin{equation}\label{eq:promom150}
\boldsymbol{\Upsilon}=\Upsilon_{ij}\,\varphi^i\otimes\varphi^j:\mathscr{Z}^2\to\mathbb{R}\quad:\Leftrightarrow\quad
(f,g)\mapsto
\boldsymbol{\Upsilon}[f,g]=\Upsilon_{ij}f^ig^j,
\end{equation}
where $\Upsilon_{ij}=\langle\varphi_i,\varphi_j\rangle$, and $i,j\in\mathbb{N}_0$.

\begin{remark}
The definition on the right-hand side of Eq.~\eqref{eq:promom150} follows from standard conventions in multi-linear algebra.
An explicit derivation is given below:
\begin{align*}
(f,g)\mapsto\boldsymbol{\Upsilon}[f,g]&=
\Upsilon_{ij}\,\varphi^i[f]\,\varphi^j[g]
=\Upsilon_{ij}\,\varphi^i[f^k\varphi_k]\,\varphi^j[g^l\varphi_l]
=\Upsilon_{ij}f^kg^l\,\varphi^i[\varphi_k]\,\varphi^j[\varphi_l]\\
&=\Upsilon_{ij}f^kg^l\delta^i_k\delta^j_l
=\Upsilon_{kl}f^kg^l
=\Upsilon_{ij}f^ig^j.
\end{align*}
The fourth and fifth equalities follow from the fact that $\varphi^i[\varphi_j]=\delta^i_j$ and that $\Upsilon_{ij}\delta^i_k\delta^j_l=\Upsilon_{kl}$ after simplification, while the last equality follows from reindexing the dummy indices $k$ and $l$ to $i$ and $j$, respectively.
The same procedure can be used to assert, for instance, that $\boldsymbol{\eta}[\mathbf{x},\mathbf{y}]=\eta_{uv}x^uy^v$ for all $\mathbf{x},\mathbf{y}\in\mathfrak{E}$, given that $\mathbf{e}^v[\mathbf{e}_u]=\delta^v_u$ by definition.
\end{remark}

Then, thanks to the orthogonality property of the random basis, the definition of $\boldsymbol{\Upsilon}$ reduces to
\begin{equation}\label{eq:promom200}
\boldsymbol{\Upsilon}[f,g]
=\Upsilon_{ij}f^ig^j
=\langle\varphi_i,\varphi_j\rangle\,f^ig^j
=\sum_{i=0}^\infty\sum_{j=0}^\infty\langle\varphi_i,\varphi_i\rangle\,\delta^i_j\, f^ig^j
=\sum_{i=0}^\infty\langle\varphi_i,\varphi_i\rangle\,f^ig^i=
\sum_{i=0}^\infty\Upsilon_{ii}f^ig^i,
\end{equation}
where the third equality follows from Eq.~\eqref{eq:spesol150}.

From Eqs.~\eqref{eq:promom100} and \eqref{eq:promom200}, we can deduce the following.
The mean of $\varrho$, $\mathbf{E}[\varrho]\in\mathbb{R}$, is simply its first component:
\begin{equation}\label{eq:promom300}
\mathbf{E}[\varrho]
=\int \varrho\,\mathrm{d}\mu
=\varrho^j\int\varphi_j\,\mathrm{d}\mu
=\varrho^j\langle\varphi_0,\varphi_j\rangle
=\Upsilon_{0j}\varrho^j=\varrho^0,
\end{equation}
while the variance of $\varrho$, $\mathrm{Var}[\varrho]\in\mathbb{R}$, is obtained by computing a specific series:
\begin{align}
\mathrm{Var}[\varrho]&=\int(\varrho-\mathbf{E}[\varrho])^2\,\mathrm{d}\mu
=\int (\varrho^j\varphi_j-\varrho^0)^2\,\mathrm{d}\mu
=\int (\varrho^j\varphi_j-\varrho^0\varphi_0)^2\,\mathrm{d}\mu\notag\\
&=\sum_{j=1}^\infty\sum_{k=1}^\infty \varrho^j\varrho^k\int\varphi_j\varphi_k\,\mathrm{d}\mu
=\sum_{j=1}^\infty\sum_{k=1}^\infty \varrho^j\varrho^k\langle\varphi_j,\varphi_k\rangle=\sum_{j=1}^\infty\Upsilon_{jj} \varrho^j\varrho^j.\label{eq:promom400}
\end{align}

\section{Numerical solution to the system of equations}\label{sec:numsol}

We employ Halley's method to numerically solve Eq.~\eqref{eq:spesol600}, thus taking advantage of its cubic convergence properties.
However, to ensure computational feasibility, we must first discretize the random function space $\mathscr{Z}$.
To this end, $\mathscr{Z}$ is $p$-discretized below, without loss of generality.

Let $\mathscr{Z}^{[P]}=\mathrm{span}(\varphi_j)_{j=0}^P$ be a finite subspace of $\mathscr{Z}$ with $P+1\in\mathbb{N}_1$ representing the dimensionality of the subspace.
If $f$ is an element of $\mathscr{Z}$, then it can be represented in $\mathscr{Z}^{[P]}$ as
\begin{equation}\label{eq:numsol100}
\sum_{j=0}^P f^j\varphi_j\equiv f^j\varphi_j.
\end{equation}

Let $Q=(R+1)(P+1)$.
If $\mathscr{Z}$ is represented with $\mathscr{Z}^{[P]}$, then Eq.~\eqref{eq:spesol600} is a system of $Q$ quadratic equations with $Q$ unknowns.
The unknowns being the two points:
\begin{equation*}
\lambda=(\lambda^0,\ldots,\lambda^i,\ldots,\lambda^P)\in\mathbb{R}^{P+1}\quad\text{and}\quad
\phi=(\phi^{10},\ldots,\phi^{ui},\ldots,\phi^{RP})\in\mathbb{R}^{R(P+1)}.
\end{equation*}

\begin{remark}
In this work, $\lambda$ has two meanings: it can represent either the random function $\lambda$ or the $(P+1)$-dimensional point $\lambda$.
Nonetheless, the intended meaning of $\lambda$ should always be clear from the context.
\end{remark}

Let $g^{ui},h^i:\mathbb{R}^Q\to\mathbb{R}$ be the associated functions of Eq.~\eqref{eq:spesol600}:
\begin{subequations}\label{eq:numsol200}
\begin{align}
g^{ui}(x)=K\indices{^u_v^i_j}\phi^{vj}-I\indices{^i_{jk}}\lambda^k\phi^{uj}&=0\label{eq:numsol200a}\\
h^i(x)=I\indices{^i_{jk}}\eta_{uv}\phi^{uj}\phi^{vk}-I\indices{^i_{00}}&=0,\label{eq:numsol200b}
\end{align}
\end{subequations}
where
\begin{equation*}
x=(\lambda,\phi)
=(\lambda^0,\ldots,\lambda^i,\ldots,\lambda^P,\phi^{10},\ldots,\phi^{ui},\ldots,\phi^{RP})\in\mathbb{R}^Q
\end{equation*}
with $u\in\{1,2,\ldots,R\}$ and $i\in\{0,1,\ldots,P\}$.

To simplify notation, the left-hand side of Eq.~\eqref{eq:numsol200} is redefined as
\begin{equation}\label{eq:numsol300}
f=g\times h:\mathbb{R}^Q\to\mathbb{R}^Q\quad:\Leftrightarrow\quad x\mapsto f(x)=(f^1(x),\ldots,f^n(x),\ldots,f^Q(x)),
\end{equation}
where $f$ is the cartesian product of $g$ and $h$ such that $f(x)=(g(x),h(x))$, $f^n:\mathbb{R}^Q\to\mathbb{R}$ is the $n$-th coordinate function of $f$ with $n\in\{1,2,\ldots,Q\}$, and the functions $g$ and $h$ are given by
\begin{align*}
g:\mathbb{R}^Q\to\mathbb{R}^{R(P+1)} &\quad:\Leftrightarrow\quad x\mapsto g(x)=(g^{10}(x),\ldots,g^{ui}(x),\ldots,g^{RP}(x))\\
h:\mathbb{R}^Q\to\mathbb{R}^{P+1} &\quad:\Leftrightarrow\quad x\mapsto h(x)=(h^0(x),\ldots,h^i(x),\ldots,h^P(x)).
\end{align*}

As a result, the system of equations to be solved hereafter is
\begin{equation}\label{eq:numsol400}
f(x)=0\in\mathbb{R}^Q\quad\text{which is equivalent to}\quad\{f^n(x)=0\}_{n=1}^Q.
\end{equation}

Linearizing Eq.~\eqref{eq:numsol400} about $x=\hat{x}$ gives
\begin{equation}\label{eq:numsol500}
f^n(\hat{x})+\frac{\partial f^n}{\partial x^m}\bigg|_{\hat{x}}(x^m-\hat{x}^m)=0,
\end{equation}
where $\hat{x}$ denotes the trial $x$, assumed to be sufficiently close to the solution, and $m\in\{1,2,\ldots,Q\}$.
(A summation is implied on the index $m$.)

Multiplying both sides of Eq.~\eqref{eq:numsol500} by
\begin{equation*}
\frac{\partial x^r}{\partial f^n}\bigg|_{f(\hat{x})}\quad\text{with $r\in\{1,2,\ldots,Q\}$},
\end{equation*}
we get
\begin{align}\label{eq:numsol600}
\frac{\partial x^r}{\partial f^n}\bigg|_{f(\hat{x})} f^n(\hat{x})+\frac{\partial x^r}{\partial f^n}\bigg|_{f(\hat{x})}\frac{\partial f^n}{\partial x^m}\bigg|_{\hat{x}} (x^m-\hat{x}^m)&=0\notag\\
\frac{\partial x^r}{\partial f^n}\bigg|_{f(\hat{x})} f^n(\hat{x})+\delta^r_m (x^m-\hat{x}^m)&=0\notag\\
\frac{\partial x^r}{\partial f^n}\bigg|_{f(\hat{x})} f^n(\hat{x})+x^r-\hat{x}^r&=0.
\end{align}
(A summation is implied on the index $n$.) The second line follows from the chain rule:
\begin{equation*}
\frac{\partial x^r}{\partial f^n}\bigg|_{f(\hat{x})}\frac{\partial f^n}{\partial x^m}\bigg|_{\hat{x}}=\frac{\partial x^r}{\partial x^m}\bigg|_{\hat{x}}=\delta^r_m.
\end{equation*}

Solving Eq.~\eqref{eq:numsol600} for $x^r-\hat{x}^r$ yields
\begin{equation}\label{eq:numsol700}
x^r-\hat{x}^r=-\frac{\partial x^r}{\partial f^n}\bigg|_{f(\hat{x})} f^n(\hat{x}).
\end{equation}

Moreover, quadratizing Eq.~\eqref{eq:numsol400} about $x=\hat{x}$ gives
\begin{equation}\label{eq:numsol800}
f^n(\hat{x})+\frac{\partial f^n}{\partial x^m}\bigg|_{\hat{x}}(x^m-\hat{x}^m)+\frac{1}{2}\frac{\partial^2f^n}{\partial x^r\partial x^m}\bigg|_{\hat{x}}(x^r-\hat{x}^r)(x^m-\hat{x}^m)=0,
\end{equation}
where $m,r\in\{1,2,\ldots,Q\}$.
(Summations are implied on the indices $m$ and $r$.)

Now let $\tilde{\mathbb{R}}^Q$ be some thickening of $\mathbb{R}^Q$ in $\mathbb{R}^{2Q}$, and let $\mathcal{H}:\tilde{\mathbb{R}}^Q\to\mathrm{L}(Q,\mathbb{R})$ be the extended Halleyian of $f$ with its $(n,m)$-th entry at $(u,v)$ defined by
\begin{equation}\label{eq:numsol900}
\mathcal{H}\indices{^n_m}(u,v)=\frac{\partial f^n}{\partial x^m}\bigg|_u+\frac{1}{2}\frac{\partial^2f^n}{\partial x^r\partial x^m}\bigg|_u v^r,
\end{equation}
where $n,m,r\in\{1,2,\ldots,Q\}$.
(A summation is implied on the index $r$.)

\begin{remark}
More specifically, in this work, the thickening of $\mathbb{R}^Q$ in $\mathbb{R}^{2Q}$ is defined as
\begin{equation*}
\tilde{\mathbb{R}}^Q=\{(u,v)\in\mathbb{R}^Q\times\mathbb{R}^Q\equiv\mathbb{R}^{2Q}:\Vert(0,v)\Vert_2<\varepsilon\}\subset\mathbb{R}^{2Q},
\end{equation*}
where $\Vert\cdot\Vert_2$ is the Euclidean norm on $\mathbb{R}^{2Q}$, and $\varepsilon>0$ is a pre-specified (small) positive number.
\end{remark}

Substituting $\hat{\mathcal{H}}\indices{^n_m}:=\mathcal{H}\indices{^n_m}(\hat{x},\,x-\hat{x})$ into Eq.~\eqref{eq:numsol800} gives
\begin{equation}\label{eq:numsol1000}
f^n(\hat{x})+\hat{\mathcal{H}}\indices{^n_m}(x^m-\hat{x}^m)=0.
\end{equation}

Thus, solving Eq.~\eqref{eq:numsol1000} for $x^m$ yields
\begin{equation}\label{eq:numsol1100}
x^m=\hat{x}^m-\hat{\mathcal{L}}\indices{^m_n}\,f^n(\hat{x}),
\end{equation}
where $\hat{\mathcal{L}}$ symbolizes the inverse of $\hat{\mathcal{H}}$, provided that $\mathcal{H}$ is non-singular at $(\hat{x},\,x-\hat{x})$; that is, now we require that $\mathcal{H}(\hat{x},\,x-\hat{x})=\hat{\mathcal{H}}\in\mathrm{GL}(Q,\mathbb{R})\subset\mathrm{L}(Q,\mathbb{R})$.

However, one problem with Eq.~\eqref{eq:numsol1100} is that $\hat{\mathcal{L}}$ depends upon $x$, the unknown in Eq.~\eqref{eq:numsol800}.
In Halley's method this problem is resolved by defining the $r$-th coordinate of $x-\hat{x}$ using Eq.~\eqref{eq:numsol700}.
Thus, by considering the function
\begin{equation*}
q:\mathbb{R}^Q\to\mathbb{R}^Q
\quad:\Leftrightarrow\quad
u\mapsto
q(u)=
(q^1(u),\ldots,q^r(u),\ldots,q^Q(u))
\quad\text{with}\quad
q^r(u)=-\frac{\partial x^r}{\partial f^s}\bigg|_{f(u)}f^s(u)\in\mathbb{R},
\end{equation*}
we can obtain the $(n,m)$-th entry of $H=\mathcal{H}(\,\cdot\,,q(\cdot)):\mathbb{R}^Q\to\mathrm{L}(Q,\mathbb{R})$ at $\hat{x}$ as follows:
\begin{align}\label{eq:numsol1200}
\hat{H}\indices{^n_m}
=H\indices{^n_m}(\hat{x})
=\mathcal{H}\indices{^n_m}(\hat{x},q(\hat{x}))
&=\frac{\partial f^n}{\partial x^m}\bigg|_{\hat{x}}-\frac{1}{2}\frac{\partial x^r}{\partial f^s}\bigg|_{f(\hat{x})}\frac{\partial^2 f^n}{\partial x^r\partial x^m}\bigg|_{\hat{x}} f^s(\hat{x})\notag\\
&=\frac{\partial f^n}{\partial x^m}\bigg|_{\hat{x}}-\frac{1}{2}\frac{\partial x^r}{\partial f^s}\bigg|_{f(\hat{x})}\frac{\partial}{\partial x^r}\left(\frac{\partial f^n}{\partial x^m}\right)\bigg|_{\hat{x}} f^s(\hat{x})\notag\\
&=\frac{\partial f^n}{\partial x^m}\bigg|_{\hat{x}}-\frac{1}{2}\frac{\partial}{\partial f^s}\left(\frac{\partial f^n}{\partial x^m}\right)\bigg|_{f(\hat{x})} f^s(\hat{x}),
\end{align}
where $H$ symbolizes the Halleyian of $f$, $s\in\{1,2,\ldots,Q\}$, and the last line follows from the chain rule.

\begin{remark}
Note that $q$ is a function that assigns to each point in $\mathbb{R}^Q$ a point in $\mathbb{O}^Q=\{v\in\mathbb{R}^Q:\Vert v\Vert_2<\varepsilon\}\subset\mathbb{R}^Q$, where $\Vert\cdot\Vert_2$ is the Euclidean norm on $\mathbb{R}^Q$, and $\varepsilon$ is the positive number pre-specified in $\tilde{\mathbb{R}}^Q=\mathbb{R}^Q\times\mathbb{O}^Q$.
\end{remark}

Hence, replacing $\hat{\mathcal{H}}$ with $\hat{H}$ in Eq.~\eqref{eq:numsol1000} results in
\begin{equation}\label{eq:numsol1140}
f^n(\hat{x})+\hat{H}\indices{^n_m}(x^m-\hat{x}^m)=0,
\end{equation}
and so, by solving for $x^m$, we get the $m$-th coordinate trial of $x$:
\begin{equation}\label{eq:numsol1150}
x^m=\hat{x}^m-\hat{L}\indices{^m_n}\,f^n(\hat{x}),
\end{equation}
where $\hat{L}$ represents the inverse of $\hat{H}$, assuming $H(\hat{x})=\hat{H}\in\mathrm{GL}(Q,\mathbb{R})$.

Supposing $f$ is locally a diffeomorphism at $\hat{x}$, let $F:\mathbb{R}^Q\to\mathrm{GL}(Q,\mathbb{R})$ and $X=F^{-1}$ be the Jacobian and the Jacobian inverse of $f$, respectively, such that their entries at $\hat{x}$ and $f(\hat{x})$ are:
\begin{equation}\label{eq:numsol1300}
\hat{F}\indices{^n_m}= 
F\indices{^n_m}(\hat{x})=
\frac{\partial f^n}{\partial x^m}\bigg|_{\hat{x}}\quad\text{and}\quad
\hat{X}\indices{^m_n}=
X\indices{^m_n}(f(\hat{x}))=
\frac{\partial x^m}{\partial f^n}\bigg|_{f(\hat{x})}.
\end{equation}

\begin{remark}
Observe that as $\hat{x}$ approaches the solution, i.e., as $f(\hat{x})\to0$, $H(\hat{x})$ tends to $F(\hat{x})$.
Indeed,
\begin{equation*}
f(\hat{x})\to 0\in\mathbb{R}^Q
\quad\Rightarrow\quad
q(\hat{x})\to 0\in\mathbb{O}^Q\subset\mathbb{R}^Q
\quad\Rightarrow\quad
H(\hat{x})\to\mathcal{H}(\hat{x},0)=F(\hat{x}).
\end{equation*}
Therefore, there is a moment at which switching from $H$ to $F$ becomes computationally advantageous. 
Nevertheless, such an approach was not implemented in this work to avoid the need for defining a stopping condition within the numerical scheme.
\end{remark}

Eq.~\eqref{eq:numsol1300} can also be expressed in matrix form as
\begin{align}\label{eq:numsol1400}
\hat{\mathbf{F}}=[\hat{F}\indices{^n_m}]=
\begin{bmatrix}
\hat{F}\indices{^{ui}_\beta} & \hat{F}\indices{^u_\alpha^i_\beta}\\[1ex]
\hat{F}\indices{^i_\beta} & \hat{F}\indices{_\alpha^i_\beta}
\end{bmatrix}\quad\text{and}\quad
\hat{\mathbf{X}}=[\hat{X}\indices{^m_n}]=
\begin{bmatrix}
\hat{X}\indices{_u^\beta_i} & \hat{X}\indices{^\beta_i}\\[1ex]
\hat{X}\indices{^\alpha_u^\beta_i} & \hat{X}\indices{^{\alpha\beta}_i}
\end{bmatrix},
\end{align}
where $n\in(ui,i)$ and $m\in(\beta,\alpha\beta)$ with $u,\alpha\in\{1,2,\ldots,R\}$ and $i,\beta\in\{0,1,\ldots,P\}$.

The entries of $\hat{\mathbf{F}}$ are determined by differentiating Eq.~\eqref{eq:numsol200} at $\hat{x}$:
\begin{subequations}\label{eq:numsol1600}
\begin{equation}\label{eq:numsol1600a}
\hat{F}\indices{^{ui}_\beta}:=
\frac{\partial g^{ui}}{\partial\lambda^\beta}\bigg|_{\hat{x}}=
-I\indices{^i_{jk}}\frac{\partial\lambda^k}{\partial\lambda^\beta}\phi^{uj}\bigg|_{\hat{x}}
=-I\indices{^i_{jk}}\delta^k_\beta\hat{\phi}^{uj}=
-I\indices{^i_{j\beta}}\hat{\phi}^{uj}=
-I\indices{^i_{\beta j}}\hat{\phi}^{uj}.
\end{equation}
\begin{align}\label{eq:numsol1600b}
\hat{F}\indices{^u_\alpha^i_\beta}:=
\frac{\partial g^{ui}}{\partial\phi^{\alpha\beta}}\bigg|_{\hat{x}}
&=K\indices{^u_v^i_j}\frac{\partial\phi^{vj}}{\partial\phi^{\alpha\beta}}\bigg|_{\hat{x}}-I\indices{^i_{jk}}\lambda^k\frac{\partial\phi^{uj}}{\partial\phi^{\alpha\beta}}\bigg|_{\hat{x}}\notag\\
&=K\indices{^u_v^i_j}\delta^v_\alpha\delta^j_\beta-I\indices{^i_{jk}}\hat{\lambda}^k\delta^u_\alpha\delta^j_\beta
=K\indices{^u_\alpha^i_\beta}-I\indices{^i_{\beta j}}\hat{\lambda}^j\delta^u_\alpha.
\end{align}
\begin{equation}\label{eq:numsol1600c}
\hat{F}\indices{^i_\beta}:=
\frac{\partial h^i}{\partial\lambda^\beta}\bigg|_{\hat{x}}=0.
\end{equation}
\begin{align}\label{eq:numsol1600d}
\hat{F}\indices{_\alpha^i_\beta}:=\frac{\partial h^i}{\partial\phi^{\alpha\beta}}\bigg|_{\hat{x}}&=
I\indices{^i_{jk}}\eta_{uv}\frac{\partial\phi^{uj}}{\partial\phi^{\alpha\beta}}\phi^{vk}\bigg|_{\hat{x}}+I\indices{^i_{jk}}\eta_{uv}\phi^{uj}\frac{\partial\phi^{vk}}{\partial\phi^{\alpha\beta}}\bigg|_{\hat{x}}\notag\\
&=I\indices{^i_{jk}}\eta_{uv}\delta^u_\alpha\delta^j_\beta\hat{\phi}^{vk}+I\indices{^i_{jk}}\eta_{uv}\hat{\phi}^{uj}\delta^v_\alpha\delta^k_\beta
=I\indices{^i_{\beta k}}\eta_{\alpha v}\hat{\phi}^{vk}+I\indices{^i_{j\beta}}\eta_{u\alpha}\hat{\phi}^{uj}\notag\\
&=I\indices{^i_{\beta j}}\eta_{\alpha u}\hat{\phi}^{uj}+I\indices{^i_{\beta j}}\eta_{\alpha u}\hat{\phi}^{uj}
=2I\indices{^i_{\beta j}}\eta_{\alpha u}\hat{\phi}^{uj}=
-2\eta_{\alpha u}\hat{F}\indices{^{ui}_\beta}.
\end{align}
\end{subequations}

The entries of $\hat{\mathbf{X}}$ are determined by taking the matrix inverse of $\hat{\mathbf{F}}$:
\begin{equation}\label{eq:numsol1650}
\hat{\mathbf{X}}=\hat{\mathbf{F}}^{-1}.
\end{equation}
Also, we remark that
\begin{equation*}
\hat{\mathbf{X}}=
\begin{bmatrix}
\hat{X}\indices{_u^\beta_i} & \hat{X}\indices{^\beta_i}\\[1ex]
\hat{X}\indices{^\alpha_u^\beta_i} & \hat{X}\indices{^{\alpha\beta}_i}
\end{bmatrix}%
:=
\begin{bmatrix}
\cfrac{\partial\lambda^\beta}{\partial g^{ui}}\bigg|_{f(\hat{x})} &
\cfrac{\partial\lambda^\beta}{\partial h^i}\bigg|_{f(\hat{x})}\\[3ex]
\cfrac{\partial\phi^{\alpha\beta}}{\partial g^{ui}}\bigg|_{f(\hat{x})} &
\cfrac{\partial\phi^{\alpha\beta}}{\partial h^i}\bigg|_{f(\hat{x})}
\end{bmatrix}.
\end{equation*}

\begin{remark}\label{rmk:numsol2000}
The correspondence between tensor components and matrix entries in Eq.~\eqref{eq:numsol1400} is as follows.
{\it For the Jacobian $\hat{\mathbf{F}}$:} $\hat{F}\indices{^{ui}_\beta}$ corresponds to the $(u+Ri,\,\beta+1)$-th entry of $\hat{\mathbf{F}}$, $\hat{F}\indices{^u_\alpha^i_\beta}$ corresponds to the $(u+Ri,\,\alpha+R\beta+P+1)$-th entry of $\hat{\mathbf{F}}$, $\hat{F}\indices{^i_\beta}$ corresponds to the $(i+R(P+1)+1,\,\beta+1)$-th entry of $\hat{\mathbf{F}}$, and $\hat{F}\indices{_\alpha^i_\beta}$ corresponds to the $(i+R(P+1)+1,\,\alpha+R\beta+P+1)$-th entry of $\hat{\mathbf{F}}$.
{\it For the Jacobian inverse $\hat{\mathbf{X}}$:} $\hat{X}\indices{_u^\beta_i}$ corresponds to the $(\beta+1,\,u+Ri)$-th entry of $\hat{\mathbf{X}}$, $\hat{X}\indices{^\beta_i}$ corresponds to the $(\beta+1,\,i+R(P+1)+1)$-th entry of $\hat{\mathbf{X}}$, $\hat{X}\indices{^\alpha_u^\beta_i}$ corresponds to the $(\alpha+R\beta+P+1,\,u+Ri)$-th entry of $\hat{\mathbf{X}}$, and $\hat{X}\indices{^{\alpha\beta}_i}$ corresponds to the $(\alpha+R\beta+P+1,\,i+R(P+1)+1)$-th entry of $\hat{\mathbf{X}}$.
\end{remark}

Another object to be determined is (from Eq.~\eqref{eq:numsol1200}):
\begin{equation}\label{eq:numsol1700}
\frac{\partial}{\partial f^s}\left(\frac{\partial f^n}{\partial x^m}\right)\bigg|_{f(\hat{x})}\equiv
\frac{\partial F\indices{^n_m}}{\partial f^s}\bigg|_{f(\hat{x})}.
\end{equation}
When expanded, this expression gives rise to eight distinct objects, whose explicit components are detailed in Appendix A.

In matrix form, the Halleyian of $f$ can be expressed as
\begin{equation}\label{eq:numsol1900}
\hat{\mathbf{H}}=[\hat{H}\indices{^n_m}]=
\begin{bmatrix}
\hat{H}\indices{^{ui}_\beta} & \hat{H}\indices{^u_\alpha^i_\beta}\\[1ex]
\hat{H}\indices{^i_\beta} & \hat{H}\indices{_\alpha^i_\beta}
\end{bmatrix},
\end{equation}
where the entries of $\hat{\mathbf{H}}$ are determined by substituting Eqs.~\eqref{eq:numsol1600} and \eqref{eq:numsol1700} into Eq.~\eqref{eq:numsol1200}.
That is:
\begin{subequations}\label{eq:numsol2000}
\begin{align}
\hat{H}\indices{^{ui}_\beta}&=
\hat{F}\indices{^{ui}_\beta}-\frac{1}{2}\frac{\partial F\indices{^{ui}_\beta}}{\partial g^{\gamma\rho}}\bigg|_{f(\hat{x})} \hat{g}^{\gamma\rho}-\frac{1}{2}\frac{\partial F\indices{^{ui}_\beta}}{\partial h^\rho}\bigg|_{f(\hat{x})} \hat{h}^\rho\notag\\
&=-I\indices{^i_{\beta j}}(\hat{\phi}^{uj}-\tfrac{1}{2}\hat{X}\indices{^u_\gamma^j_\rho}\hat{g}^{\gamma\rho}-\tfrac{1}{2}\hat{X}\indices{^{uj}_\rho}\hat{h}^\rho).
\end{align}
\begin{align}
\hat{H}\indices{^u_\alpha^i_\beta}&=
\hat{F}\indices{^u_\alpha^i_\beta}-\frac{1}{2}\frac{\partial F\indices{^u_\alpha^i_\beta}}{\partial g^{\gamma\rho}}\bigg|_{f(\hat{x})} \hat{g}^{\gamma\rho}-\frac{1}{2}\frac{\partial F\indices{^u_\alpha^i_\beta}}{\partial h^\rho}\bigg|_{f(\hat{x})} \hat{h}^\rho\notag\\
&=K\indices{^u_\alpha^i_\beta}-I\indices{^i_{\beta j}}(\hat{\lambda}^j-\tfrac{1}{2}\hat{X}\indices{_\gamma^j_\rho}\hat{g}^{\gamma\rho}-\tfrac{1}{2}\hat{X}\indices{^j_\rho}\hat{h}^\rho)\,\delta^u_\alpha.
\end{align}
\begin{equation}
\hat{H}\indices{^i_\beta}=
\hat{F}\indices{^i_\beta}-\frac{1}{2}\frac{\partial F\indices{^i_\beta}}{\partial g^{\gamma\rho}}\bigg|_{f(\hat{x})}\hat{g}^{\gamma\rho}-\frac{1}{2}\frac{\partial F\indices{^i_\beta}}{\partial h^\rho}\bigg|_{f(\hat{x})}\hat{h}^\rho=0.
\end{equation}
\begin{align}
\hat{H}\indices{_\alpha^i_\beta}&=
\hat{F}\indices{_\alpha^i_\beta}-\frac{1}{2}\frac{\partial F\indices{_\alpha^i_\beta}}{\partial g^{\gamma\rho}}\bigg|_{f(\hat{x})}\hat{g}^{\gamma\rho}-\frac{1}{2}\frac{\partial F\indices{_\alpha^i_\beta}}{\partial h^\rho}\bigg|_{f(\hat{x})}\hat{h}^\rho\notag\\
&=2\eta_{\alpha u}I\indices{^i_{\beta j}}(\hat{\phi}^{uj}-\tfrac{1}{2}\hat{X}\indices{^u_\gamma^j_\rho}\hat{g}^{\gamma\rho}-\tfrac{1}{2}\hat{X}\indices{^{uj}_\rho}\hat{h}^\rho)
=-2\eta_{\alpha u}\hat{H}\indices{^{ui}_\beta}.
\end{align}
\end{subequations}

These expressions can be further shortened if these objects are defined:
\begin{subequations}\label{eq:numsol2100}
\begin{gather}
\bar{\lambda}^j=\hat{\lambda}^j-\tfrac{1}{2}\hat{X}\indices{_\gamma^j_\rho}\hat{g}^{\gamma\rho}-\tfrac{1}{2}\hat{X}\indices{^j_\rho}\hat{h}^\rho\\
\bar{\phi}^{uj}=\hat{\phi}^{uj}-\tfrac{1}{2}\hat{X}\indices{^u_\gamma^j_\rho}\hat{g}^{\gamma\rho}-\tfrac{1}{2}\hat{X}\indices{^{uj}_\rho}\hat{h}^\rho.
\end{gather}
\end{subequations}

Hence, if Eq.~\eqref{eq:numsol2100} is replaced into Eq.~\eqref{eq:numsol2000}, the entries of $\hat{\mathbf{H}}$ become
\begin{subequations}\label{eq:numsol2200}
\begin{gather}
\hat{H}\indices{^{ui}_\beta}=-I\indices{^i_{\beta j}}\bar{\phi}^{uj}\\
\hat{H}\indices{^u_\alpha^i_\beta}=K\indices{^u_\alpha^i_\beta}-I\indices{^i_{\beta j}}\bar{\lambda}^j\delta^u_\alpha\\
\hat{H}\indices{^i_\beta}=0\\
\hat{H}\indices{_\alpha^i_\beta}=-2\eta_{\alpha u}\hat{H}\indices{^{ui}_\beta}.
\end{gather}
\end{subequations}

Some useful remarks are in order.

\begin{remark}\label{rmk:numsol3000}
The correspondence between tensor components and matrix entries in Eq.~\eqref{eq:numsol2200} is the same as that of $\hat{F}$ (i.e., Remark \ref{rmk:numsol2000} with $\hat{F}$ substituted for $\hat{H}$).
\end{remark}

\begin{remark}\label{rmk:numsol4000}
From a computational standpoint, Eq.~\eqref{eq:numsol2100} is best written in matrix form as
\begin{equation}\label{eq:numsol2300}
\bar{\mathbf{x}}=
\big\{\bar{x}^n\big\}=
\begin{Bmatrix}
\bar{\lambda}^j\\[1ex]
\bar{\phi}^{uj}
\end{Bmatrix}=
\hat{\mathbf{x}}-\tfrac{1}{2}\hat{\mathbf{X}}\hat{\mathbf{f}},
\end{equation}
where $u\in\{1,2,\ldots,R\}$, $j\in\{0,1,\ldots,P\}$,
\begin{equation}\label{eq:numsol2400}
\hat{\mathbf{x}}=
\big\{\hat{x}^n\big\}=
\begin{Bmatrix}
\hat{\lambda}^j\\[1ex]
\hat{\phi}^{uj}
\end{Bmatrix}\quad\text{and}\quad
\hat{\mathbf{f}}=
\big\{\hat{f}^n\big\}=
\begin{Bmatrix}
\hat{g}^{uj}\\[1ex]
\hat{h}^j
\end{Bmatrix}.
\end{equation}
\end{remark}

\begin{remark}\label{rmk:numsol5000}
The correspondence between vector components and vector entries in Eq.~\eqref{eq:numsol2400} is as follows: $\hat{\lambda}^j$ corresponds to the $(j+1)$-th entry of $\hat{\mathbf{x}}$, $\hat{\phi}^{uj}$ corresponds to the $(u+Rj+P+1)$-th entry of $\hat{\mathbf{x}}$, $\hat{g}^{uj}$ corresponds to the $(u+Rj)$-th entry of $\hat{\mathbf{f}}$, and $\hat{h}^j$ corresponds to the $(j+R(P+1)+1)$-th entry of $\hat{\mathbf{f}}$.
\end{remark}

\begin{remark}\label{rmk:numsol6000}
It is noted that when Eq.~\eqref{eq:numsol2200} is compared to Eq.~\eqref{eq:numsol1600}, one gets:
\begin{equation}\label{eq:numsol2500}
\hat{H}\indices{^n_m}=H\indices{^n_m}(\hat{x})=F\indices{^n_m}(\bar{x}).
\end{equation}
In other terms, $\hat{\mathbf{H}}=\bar{\mathbf{F}}$.
\end{remark}

\section{Error analysis of proposed method}\label{sec:errana}

In this section, all indices range over $\{1,2,\ldots,Q\}$, unless otherwise noted.
If any indices are repeated within a term, summation over the repeated index is implied.

Let $y$ be a root of $f$, as defined by Eq.~\eqref{eq:numsol400}.
Then $f(y)=0$.
If $f$ is three times continuously differentiable at $\hat{x}$, and $y$ is sufficiently close to $\hat{x}$, then Taylor's theorem implies that there exist points $\xi$ and $\eta$ within the region bounded by $\hat{x}$ and $y$ (two diametrically opposite points defining the opposite vertices of a hyperrectangle), such that:
\begin{subequations}\label{eq:errana1000}
\begin{multline}\label{eq:errana1000a}
f^n(\hat{x})
+\frac{\partial f^n}{\partial x^m}\bigg|_{\hat{x}}(y^m-\hat{x}^m)
+\frac{1}{2}\frac{\partial^2 f^n}{\partial x^r\partial x^m}\bigg|_{\hat{x}}(y^r-\hat{x}^r)(y^m-\hat{x}^m)\\
+\frac{1}{6}\frac{\partial^3 f^n}{\partial x^a\partial x^b\partial x^m}\bigg|_\xi
(y^a-\hat{x}^a)(y^b-\hat{x}^b)(y^m-\hat{x}^m)=0
\end{multline}
and
\begin{equation}\label{eq:errana1000b}
f^s(\hat{x})
+\frac{\partial f^s}{\partial x^m}\bigg|_{\hat{x}}(y^m-\hat{x}^m)
+\frac{1}{2}\frac{\partial^2 f^s}{\partial x^a\partial x^b}\bigg|_\eta(y^a-\hat{x}^a)(y^b-\hat{x}^b)=0.
\end{equation}
\end{subequations}

Multiplying both sides of Eq.~\eqref{eq:errana1000b} by
\begin{equation*}
\frac{\partial x^r}{\partial f^s}\bigg|_{f(\hat{x})}
\end{equation*}
and simplifying produces
\begin{equation}\label{eq:errana2000}
\frac{\partial x^r}{\partial f^s}\bigg|_{f(\hat{x})}f^s(\hat{x})
+y^r-\hat{x}^r
+\frac{1}{2}\frac{\partial x^r}{\partial f^s}\bigg|_{f(\hat{x})}\frac{\partial^2 f^s}{\partial x^a\partial x^b}\bigg|_\eta(y^a-\hat{x}^a)(y^b-\hat{x}^b)=0.
\end{equation}
Thus, solving for $y^r-\hat{x}^r$ in Eq.~\eqref{eq:errana2000} yields
\begin{equation}\label{eq:errana3000}
y^r-\hat{x}^r=
-\frac{\partial x^r}{\partial f^s}\bigg|_{f(\hat{x})}f^s(\hat{x})
-\frac{1}{2}\frac{\partial x^r}{\partial f^s}\bigg|_{f(\hat{x})}\frac{\partial^2 f^s}{\partial x^a\partial x^b}\bigg|_\eta(y^a-\hat{x}^a)(y^b-\hat{x}^b).
\end{equation}

Moreover, Eq.~\eqref{eq:errana1000a} can also be expressed as
\begin{equation}\label{eq:errana4000}
f^n(\hat{x})
+\bigg[\frac{\partial f^n}{\partial x^m}\bigg|_{\hat{x}}
+\frac{1}{2}\frac{\partial^2 f^n}{\partial x^r\partial x^m}\bigg|_{\hat{x}}(y^r-\hat{x}^r)\bigg](y^m-\hat{x}^m)
+\frac{1}{6}\frac{\partial^3 f^n}{\partial x^a\partial x^b\partial x^m}\bigg|_\xi
(y^a-\hat{x}^a)(y^b-\hat{x}^b)(y^m-\hat{x}^m)=0.
\end{equation}
Inserting Eq.~\eqref{eq:errana3000} into the brackets of Eq.~\eqref{eq:errana4000} and rearranging terms gives
\begin{multline}\label{eq:errana5000}
f^n(\hat{x})
+\bigg[\frac{\partial f^n}{\partial x^m}\bigg|_{\hat{x}}
-\frac{1}{2}\frac{\partial x^r}{\partial f^s}\bigg|_{f(\hat{x})}\frac{\partial^2 f^n}{\partial x^r\partial x^m}\bigg|_{\hat{x}}f^s(\hat{x})\bigg](y^m-\hat{x}^m)\\
-\frac{1}{4}\frac{\partial x^r}{\partial f^s}\bigg|_{f(\hat{x})}\frac{\partial^2 f^n}{\partial x^r\partial x^m}\bigg|_{\hat{x}}\frac{\partial^2 f^s}{\partial x^a\partial x^b}\bigg|_\eta(y^a-\hat{x}^a)(y^b-\hat{x}^b)(y^m-\hat{x}^m)\\
+\frac{1}{6}\frac{\partial^3 f^n}{\partial x^a\partial x^b\partial x^m}\bigg|_\xi
(y^a-\hat{x}^a)(y^b-\hat{x}^b)(y^m-\hat{x}^m)=0.
\end{multline}
Reindexing, simplifying and rearranging terms (one more time) yields
\begin{multline}\label{eq:errana6000}
f^n(\hat{x})
+\bigg[\frac{\partial f^n}{\partial x^m}\bigg|_{\hat{x}}
-\frac{1}{2}\frac{\partial}{\partial f^s}\bigg(\frac{\partial f^n}{\partial x^m}\bigg)\bigg|_{f(\hat{x})}f^s(\hat{x})\bigg](y^m-\hat{x}^m)=\\
\bigg[\frac{1}{4}\frac{\partial^2 f^s}{\partial x^a\partial x^b}\bigg|_\eta\frac{\partial}{\partial f^s}\bigg(\frac{\partial f^n}{\partial x^c}\bigg)\bigg|_{f(\hat{x})}
-\frac{1}{6}\frac{\partial^3 f^n}{\partial x^a\partial x^b\partial x^c}\bigg|_\xi\bigg]
(y^a-\hat{x}^a)(y^b-\hat{x}^b)(y^c-\hat{x}^c),
\end{multline}
where the left-hand side of the equation is nothing but: $f^n(\hat{x})+\hat{H}\indices{^n_m}(y^m-\hat{x}^m)$.

Subtracting Eq.~\eqref{eq:numsol1140} from Eq.~\eqref{eq:errana6000} and then solving for $y^m-x^m$ gives
\begin{equation}\label{eq:errana7000}
y^m-x^m=
\hat{L}\indices{^m_n}\bigg[\frac{1}{4}\frac{\partial^2 f^s}{\partial x^a\partial x^b}\bigg|_\eta\frac{\partial}{\partial f^s}\bigg(\frac{\partial f^n}{\partial x^c}\bigg)\bigg|_{f(\hat{x})}
-\frac{1}{6}\frac{\partial^3 f^n}{\partial x^a\partial x^b\partial x^c}\bigg|_\xi\bigg]
(y^a-\hat{x}^a)(y^b-\hat{x}^b)(y^c-\hat{x}^c).
\end{equation}

Thus, as $\hat{x}\to y$, also $\xi,\eta\to y$ and
\begin{multline}\label{eq:errana8000}
\hat{L}\indices{^m_n}\bigg[\frac{1}{4}\frac{\partial^2 f^s}{\partial x^a\partial x^b}\bigg|_\eta\frac{\partial}{\partial f^s}\bigg(\frac{\partial f^n}{\partial x^c}\bigg)\bigg|_{f(\hat{x})}
-\frac{1}{6}\frac{\partial^3 f^n}{\partial x^a\partial x^b\partial x^c}\bigg|_\xi\bigg]\\
\to\check{E}\indices{^m_{abc}}:=
\check{L}\indices{^m_n}\bigg[\frac{1}{4}\frac{\partial^2 f^s}{\partial x^a\partial x^b}\bigg|_y\frac{\partial}{\partial f^s}\bigg(\frac{\partial f^n}{\partial x^c}\bigg)\bigg|_{f(y)}
-\frac{1}{6}\frac{\partial^3 f^n}{\partial x^a\partial x^b\partial x^c}\bigg|_y\bigg],
\end{multline}
where $\check{L}$ is the inverse of $H$ at $y$ (provided is non-singular).

Substituting the converging constant $\check{E}\indices{^m_{abc}}$ into Eq.~\eqref{eq:errana7000} gives
\begin{equation}\label{eq:errana9000}
y^m-x^m\sim\check{E}\indices{^m_{abc}}(y^a-\hat{x}^a)(y^b-\hat{x}^b)(y^c-\hat{x}^c).
\end{equation}

Therefore, for the eigenvalue problem at hand, the following conclusions are noted.
First, since the eigenvalue problem is inherently quadratic, it is clear that
\begin{equation}\label{eq:errana11000}
\frac{\partial^3 f^n}{\partial x^a\partial x^b\partial x^c}\bigg|_y=0.
\end{equation}
Second, the expression:
\begin{equation}\label{eq:errana12000}
\frac{\partial}{\partial f^s}\bigg(\frac{\partial f^n}{\partial x^c}\bigg)\bigg|_{f(y)}\equiv
\frac{\partial F\indices{^n_c}}{\partial f^s}\bigg|_{f(y)}
\end{equation}
is the same as the one presented earlier in Eq.~\eqref{eq:numsol1700}.
Third, the expansion of
\begin{equation}\label{eq:errana13000}
\frac{\partial^2 f^s}{\partial x^a\partial x^b}\bigg|_y\equiv\frac{\partial F\indices{^s_b}}{\partial x^a}\bigg|_y
\end{equation}
yields eight objects, whose components are formally derived in Appendix B.

Thus, upon expanding
\begin{equation}\label{eq:errana15000}
\check{D}\indices{^n_{abc}}:=\frac{\partial^2 f^s}{\partial x^a\partial x^b}\bigg|_y\frac{\partial}{\partial f^s}\bigg(\frac{\partial f^n}{\partial x^c}\bigg)\bigg|_{f(y)},
\end{equation}
it gives rise to 16 distinct objects.
The derivation of each individual component is provided in Appendix C.

Moreover, let $\Delta x^n=y^n-\hat{x}^n$. 
Then, from Eqs.~\eqref{eq:errana8000}, \eqref{eq:errana9000} and \eqref{eq:errana15000}, we obtain:
\begin{align}\label{eq:errana17000}
y^m-x^m &\sim\check{E}\indices{^m_{abc}}\Delta x^a\Delta x^b\Delta x^c=\tfrac{1}{4}\check{L}\indices{^m_n}\check{D}\indices{^n_{abc}}\Delta x^a\Delta x^b\Delta x^c\notag\\
&=\tfrac{1}{4}\sum_{n\in(vj,j)}\sum_{a\in(\rho,\gamma\rho)}\sum_{b\in(\beta,\alpha\beta)}\sum_{c\in(\sigma,\omega\sigma)}\check{L}\indices{^m_n}\check{D}\indices{^n_{abc}}\Delta x^a\Delta x^b\Delta x^c.
\end{align}

Thus, substituting Eq.~\eqref{eq:errana16000} into Eq.~\eqref{eq:errana17000}, recalling that $\boldsymbol{\tilde{I}}$ is symmetric with respect to its subscripts ($I\indices{^i_{jk}}=I\indices{^i_{kj}}$), reindexing, and simplifying, yields
\begin{multline}\label{eq:errana18000}
y^m-x^m\sim
\check{A}\indices{^m_{\rho\beta(\omega\sigma)}}\Delta\lambda^\rho\Delta\lambda^\beta\Delta\phi^{\omega\sigma}
+\check{A}\indices{^m_{\rho(\alpha\beta)(\omega\sigma)}}\Delta\lambda^\rho\Delta\phi^{\alpha\beta}\Delta\phi^{\omega\sigma}\\
+\check{A}\indices{^m_{(\gamma\rho)(\alpha\beta)(\omega\sigma)}}\Delta\phi^{\gamma\rho}\Delta\phi^{\alpha\beta}\Delta\phi^{\omega\sigma},
\end{multline}
where
\begin{subequations}
\begin{equation}
\check{A}\indices{^m_{\rho\beta(\omega\sigma)}}=\tfrac{1}{2}I\indices{^i_{\sigma\rho}}I\indices{^j_{\beta k}}\check{L}\indices{^m_{(vj)}}\check{X}\indices{^v_\omega^k_i},
\end{equation}
\begin{equation}
\check{A}\indices{^m_{\rho(\alpha\beta)(\omega\sigma)}}=-\tfrac{1}{2}I\indices{^i_{\beta\rho}}I\indices{^j_{\sigma k}}(2\eta_{\omega w}\check{L}\indices{^m_j}\check{X}\indices{^w_\alpha^k_i}-\check{L}\indices{^m_{(\omega j)}}\check{X}\indices{_\alpha^k_i})
-\tfrac{1}{2}I\indices{^i_{\beta\sigma}}I\indices{^j_{\rho k}}\eta_{\alpha\omega}\check{L}\indices{^m_{(vj)}}\check{X}\indices{^{vk}_i},
\end{equation}
\begin{equation}
\check{A}\indices{^m_{(\gamma\rho)(\alpha\beta)(\omega\sigma)}}=\tfrac{1}{2}I\indices{^i_{\beta\rho}}I\indices{^j_{\sigma k}}\eta_{\alpha\gamma}(2\eta_{\omega w}\check{L}\indices{^m_j}\check{X}\indices{^{wk}_i}-\check{L}\indices{^m_{(\omega j)}}\check{X}\indices{^k_i}).
\end{equation}
\end{subequations}

Consequently, the error associated with the $m$-th coordinate trial of $x$ is a linear combination of the cube of the previous errors: $\Delta\lambda^\rho\Delta\lambda^\beta\Delta\phi^{\omega\sigma}$, $\Delta\lambda^\rho\Delta\phi^{\alpha\beta}\Delta\phi^{\omega\sigma}$ and $\Delta\phi^{\gamma\rho}\Delta\phi^{\alpha\beta}\Delta\phi^{\omega\sigma}$, as computed by Eq.~\eqref{eq:errana18000}, assuming that $H$ at $y$ is non-singular.
It is noted that the error is not a linear combination of $\Delta\lambda^\rho\Delta\lambda^\beta\Delta\lambda^\sigma$, which implies that once the $\phi^{ui}$'s are known, the corresponding $\lambda^i$'s will immediately be known in the next iteration.
This underscores the significance of using Halley's method as the root-finding algorithm for Eq.~\eqref{eq:spesol600}, as it can significantly accelerate the solution process when the eigenvector components, $\phi^{ui}$'s, are nearly known.

\section{Computational cost of proposed method}\label{sec:comcos}

In this section, we present the computational cost of Halley's method by evaluating its time complexity per iteration using Eq.~\eqref{eq:numsol1150}.
Given $\boldsymbol{K}$, $\boldsymbol{I}$ and $\boldsymbol{\eta}$, the goal of Halley's scheme is to calculate the next trial solution $\mathbf{x}$, based on the current trial solution $\hat{\mathbf{x}}$:
\begin{equation}\label{eq:comcos500}
x^m=\hat{x}^m+(-1)\sum_{n=1}^Q\hat{L}\indices{^m_n}\hat{f}^n
\quad\Leftrightarrow\quad
\mathbf{x}=\hat{\mathbf{x}}+(-1)\hat{\mathbf{L}}\hat{\mathbf{f}},
\end{equation}
where $\hat{f}^n=f^n(\hat{x})$, and $m\in\{1,2,\ldots,Q\}$.
The equation on the right-hand side is the matrix representation of the one on the left-hand side.
Hence, $\hat{\mathbf{x}}$ and $\hat{\mathbf{f}}$ are vectors of size $Q$, and $\hat{\mathbf{L}}$ is a square matrix of size $Q$.

Now, if $\hat{\mathbf{x}}$, $\hat{\mathbf{L}}$ and $\hat{\mathbf{f}}$ are known, the time complexity of $\mathbf{x}$ is $\mathrm{tc}'(\mathbf{x})=\mathcal{O}(Q^2)=\mathcal{O}(R^2P^2)$.
This follows from the fact that $Q=(R+1)(P+1)$, and that the total number of operations required to calculate $\mathbf{x}$ is 
\begin{equation}
\mathrm{no}'(\mathbf{x})=Q(2Q+1)=2Q^2+Q.
\end{equation}

The primes in $\mathrm{tc}'(\mathbf{x})$ and $\mathrm{no}'(\mathbf{x})$ indicate that the corresponding expressions are valid under the assumption that $\hat{\mathbf{x}}$, $\hat{\mathbf{L}}$ and $\hat{\mathbf{f}}$ have been precomputed.
To eliminate this assumption, and compute the time complexity of $\mathbf{x}$ appropriately, it is necessary to compute firstly the matrices $\hat{\mathbf{F}}$, $\hat{\mathbf{X}}$, $\hat{\mathbf{H}}$ and $\hat{\mathbf{L}}$ and the vector $\hat{\mathbf{f}}$.
This is done below.

Starting with the time complexities of $\hat{\mathbf{F}}$ and $\hat{\mathbf{X}}$ ($=\hat{\mathbf{F}}^{-1}$), we refer to Eq.~\eqref{eq:numsol1600} to observe that:
\begin{subequations}\label{eq:comcos1000}
\begin{equation}\label{eq:comcos1000a}
\hat{F}\indices{^{ui}_\beta}=
(-1)\sum_{j=0}^P I\indices{^i_{\beta j}}\hat{\phi}^{uj}
\quad\text{implies}\quad
2(P+1)\text{ operations}.
\end{equation}
\begin{equation}\label{eq:comcos1000b}
\hat{F}\indices{^u_\alpha^i_\beta}=
K\indices{^u_\alpha^i_\beta}+(-1)\delta^u_\alpha\sum_{j=0}^PI\indices{^i_{\beta j}}\hat{\lambda}^j
\quad\text{implies}\quad
2(P+2)\text{ operations}.
\end{equation}
\begin{equation}\label{eq:comcos1000c}
\hat{F}\indices{^i_\beta}=0
\quad\text{implies}\quad
0\text{ operations}.
\end{equation}
\begin{equation}\label{eq:comcos1000d}
\hat{F}\indices{_\alpha^i_\beta}=
(-2)\sum_{u=1}^R\eta_{\alpha u}\hat{F}\indices{^{ui}_\beta}
\quad\text{implies}\quad
2R\text{ operations}.
\end{equation}
\end{subequations}

\begin{remark}
To determine the number of operations in Eq.~\eqref{eq:comcos1000a}, we note that the summation over $j=\{0,1,\ldots,P\}$ involves $P+1$ multiplications and $P$ additions.
After the summation, the result is multiplied by $-1$, adding one more multiplication. 
Therefore, a total of $2(P+1)$ arithmetic operations are required to compute $\hat{F}\indices{^{ui}_\beta}$.
The same procedure was used to obtain the number of operations for the remaining components of $\hat{\mathbf{F}}$, and it is likewise used to compute those for $\hat{\mathbf{H}}$ and $\hat{\mathbf{f}}$.
\end{remark}

From Remark \ref{rmk:numsol2000}, we observe that in $\hat{\mathbf{F}}$: $\hat{F}\indices{^{ui}_\beta}$ represents $R(P+1)^2$ entries, $\hat{F}\indices{^u_\alpha^i_\beta}$ represents $R^2(P+1)^2$ entries, $\hat{F}\indices{^i_\beta}$ represents $(P+1)^2$ entries, and $\hat{F}\indices{_\alpha^i_\beta}$ represents $R(P+1)^2$ entries.
Taking this observation into account, we obtain the total number of operations required to construct the matrix $\hat{\mathbf{F}}$:
\begin{align}
\mathrm{no}(\hat{\mathbf{F}})&=R(P+1)^2\cdot 2(P+1) + R^2(P+1)^2\cdot 2(P+2) + (P+1)^2\cdot 0 + R(P+1)^2\cdot 2R\notag\\
&=2R(P+1)^2\big[(R+1)(P+1)+2R\big].
\end{align}
Hence, the time complexity of matrix $\hat{\mathbf{F}}$ is $\mathrm{tc}(\hat{\mathbf{F}})=\mathcal{O}(R^2P^3)$.

The time complexity of inverting matrix $\hat{\mathbf{F}}$ using Gauss-Jordan elimination is $\mathrm{tc}(\hat{\mathbf{F}})=\mathcal{O}(Q^3)$, which means that the time complexity of matrix $\hat{\mathbf{X}}$ is $\mathrm{tc}(\hat{\mathbf{X}})=\mathcal{O}(R^3P^3)$.

Next, we compute the time complexities of $\hat{\mathbf{H}}$ and $\hat{\mathbf{L}}$ ($=\hat{\mathbf{H}}^{-1}$).
From Eq.~\eqref{eq:numsol2000},
\begin{subequations}\label{comcos2000}
\begin{multline}
\hat{H}\indices{^{ui}_\beta}=
(-1)\sum_{j=0}^P I\indices{^i_{\beta j}}\!\bigg[\hat{\phi}^{uj}+(-\tfrac{1}{2})\bigg(\sum_{\gamma=1}^R\sum_{\rho=0}^P\hat{X}\indices{^u_\gamma^j_\rho}\hat{g}^{\gamma\rho}+\sum_{\rho=0}^P\hat{X}\indices{^{uj}_\rho}\hat{h}^\rho\bigg)\bigg]\\
\quad\text{implies}\quad
2(R+2)(P+1)+1\text{ operations}.
\end{multline}
\begin{multline}
\hat{H}\indices{^u_\alpha^i_\beta}=
K\indices{^u_\alpha^i_\beta}+(-1)\delta^u_\alpha\sum_{j=0}^P I\indices{^i_{\beta j}}\!\bigg[\hat{\lambda}^j+(-\tfrac{1}{2})\bigg(\sum_{\gamma=1}^R\sum_{\rho=0}^P\hat{X}\indices{_\gamma^j_\rho}\hat{g}^{\gamma\rho}+\sum_{\rho=0}^P\hat{X}\indices{^j_\rho}\hat{h}^\rho\bigg)\bigg]\\
\quad\text{implies}\quad
2(R+2)(P+1)+3\text{ operations}.
\end{multline}
\begin{equation}
\hat{H}\indices{^i_\beta}=0
\quad\text{implies}\quad
0\text{ operations}.
\end{equation}
\begin{equation}
\hat{H}\indices{_\alpha^i_\beta}=
(-2)\sum_{u=1}^R\eta_{\alpha u}\hat{H}\indices{^{ui}_\beta}
\quad\text{implies}\quad
2R\text{ operations}.
\end{equation}
\end{subequations}

Thus, the total number of operations required to construct the matrix $\hat{\mathbf{H}}$, after recognizing that its components (see Remark \ref{rmk:numsol3000}) share the same structure as those of $\hat{\mathbf{F}}$, is:
\begin{equation}
\mathrm{no}(\hat{\mathbf{H}})=R(P+1)^2\big[2(R+2)(R+1)(P+1)+5R+1\big],
\end{equation}
meaning that the time complexity of $\hat{\mathbf{H}}$ is $\mathrm{tc}(\hat{\mathbf{H}})=\mathcal{O}(R^3P^3)$ and that of $\hat{\mathbf{L}}$ (using again Gauss-Jordan elimination as the method for matrix inversion) is $\mathrm{tc}(\hat{\mathbf{L}})=\mathcal{O}(R^3P^3)$.

Lastly, we compute the time complexity of $\hat{\mathbf{f}}$.
From Eq.~\eqref{eq:numsol200}, we note that:
\begin{subequations}
\begin{equation}
\hat{g}^{ui}=\sum_{v=1}^R\sum_{j=0}^P K\indices{^u_v^i_j}\hat{\phi}^{vj}+(-1)\sum_{j=0}^P\sum_{k=0}^P I\indices{^i_{jk}}\hat{\lambda}^k\hat{\phi}^{uj}
\quad\text{implies}\quad
2(P+1)(R+P+1)\text{ operations}.
\end{equation}
\begin{equation}
\hat{h}^i=(-1)I\indices{^i_{00}}+\sum_{u=1}^R\sum_{v=1}^R\sum_{j=0}^P\sum_{k=0}^P I\indices{^i_{jk}}\eta_{uv}\hat{\phi}^{uj}\hat{\phi}^{vk}
\quad\text{implies}\quad
2R^2(P+1)^2+1\text{ operations}.
\end{equation}
\end{subequations}

From Remark \ref{rmk:numsol5000}, we see that in $\hat{\mathbf{f}}$: $\hat{g}^{ui}$ represents $R(P+1)$ entries, and $\hat{h}^i$ represents $P+1$ entries.
Hence, the total number of operations required to construct the vector $\hat{\mathbf{f}}$ is:
\begin{align}
\mathrm{no}(\hat{\mathbf{f}})&=R(P+1)\cdot 2(P+1)(R+P+1) + (P+1)\cdot(2R^2(P+1)^2+1)\notag\\
&=(P+1)\big[2R(R+1)(P+1)^2+2R^2(P+1)+1\big],
\end{align}
which means that the time complexity of $\hat{\mathbf{f}}$ is $\mathrm{tc}(\hat{\mathbf{f}})=\mathcal{O}(R^2P^3)$.

As a result, in Halley's method, the time complexity of calculating $\mathbf{x}$ (the next trial solution) is:
\begin{align}
\mathrm{tc}\text{-Halley}(\mathbf{x})&=\mathrm{tc}(\hat{\mathbf{F}})+\mathrm{tc}(\hat{\mathbf{X}})+\mathrm{tc}(\hat{\mathbf{H}})+\mathrm{tc}(\hat{\mathbf{L}})+\mathrm{tc}(\hat{\mathbf{f}})+\mathrm{tc}'(\mathbf{x})\notag\\
&=\mathcal{O}(R^2P^3)+\mathcal{O}(R^3P^3)+\mathcal{O}(R^3P^3)+\mathcal{O}(R^3P^3)+\mathcal{O}(R^2P^3)+\mathcal{O}(R^2P^2)\notag\\
&=\mathcal{O}(R^3P^3).
\end{align}
Interestingly, this is the same time complexity produced by Newton's method (with $\hat{\mathbf{L}}=\hat{\mathbf{X}}$ in Eq.~\eqref{eq:comcos500}):
\begin{align}
\mathrm{tc}\text{-Newton}(\mathbf{x})&=\mathrm{tc}(\hat{\mathbf{F}})+\mathrm{tc}(\hat{\mathbf{X}})+\mathrm{tc}(\hat{\mathbf{f}})+\mathrm{tc}'(\mathbf{x})\notag\\
&=\mathcal{O}(R^2P^3)+\mathcal{O}(R^3P^3)+\mathcal{O}(R^2P^3)+\mathcal{O}(R^2P^2)\notag\\
&=\mathcal{O}(R^3P^3).
\end{align}

However, it is important to add that the primary computational bottleneck in Newton's method solely arises from computing $\hat{\mathbf{X}}$, whereas in Halley's method, it stems from computing $\hat{\mathbf{X}}$, $\hat{\mathbf{H}}$ and $\hat{\mathbf{L}}$.
As a result, a single iteration of Halley's method is approximately three times more computationally expensive than one iteration of Newton's method.
Nevertheless, as demonstrated in the previous section, Halley's method offers a significant advantage in convergence rate, and the potential to further accelerate the solution process when the eigenvector components are roughly known.

\section{Illustrative example}\label{sec:illexa}

Consider the problem of solving the 3-by-3 system of stochastic linear ODEs
\begin{equation}\label{eq:illexa1000}
\dot{\mathbf{x}}=\mathbf{K}\mathbf{x}
\qquad:\Leftrightarrow\qquad
\begin{bmatrix}
\dot{x}\\
\dot{y}\\
\dot{z}
\end{bmatrix}=
\begin{bmatrix}
2b & -b & 0\\
-b & b+c & -c\\
0 & -c & c
\end{bmatrix}\!
\begin{bmatrix}
x\\
y\\
z
\end{bmatrix},
\end{equation}
where $b,c:\mathfrak{Z}\to\mathbb{R}$ are random variables given by 
\begin{equation*}
b(\zeta)=5\zeta^1+20\quad\text{and}\quad c(\zeta)=10\zeta^2+15
\end{equation*}
(with the superscripts not denoting exponentiation), and $\zeta=(\zeta^1,\zeta^2)\sim\mathrm{Uniform}^{\otimes 2}$ is a two-tuple random variable uniformly distributed in $\mathfrak{Z}=[-1,1]^2$, where the probability density function, $f:\mathfrak{Z}\to\mathbb{R}_0^+$, is naturally given by $f(\zeta)=\tfrac{1}{4}$.
As a result, the measure on $(\mathfrak{Z},\boldsymbol{\mathfrak{Z}})$ is given by:
\begin{equation*}
\mu(\mathfrak{A})=\int_\mathfrak{A} f(\zeta)\,\mathrm{d}\zeta=\tfrac{1}{4}\int_\mathfrak{A}\mathrm{d}\zeta,\quad\text{for all $\mathfrak{A}\in\boldsymbol{\mathfrak{Z}}$}.
\end{equation*}

Since one solution of Eq.~\eqref{eq:illexa1000} is $\mathbf{x}=\boldsymbol{\upphi}e^{\lambda t}$, the problem reduces to solving Eq.~\eqref{eq:prosta100}, after noticing that $\mathbf{K}$ is a random, real symmetric matrix of finite size ($R=3$).
For convenience, $\mathbf{K}:\mathfrak{Z}\to\mathrm{Sym}(3,\mathbb{R})$ is taken hereafter as
\begin{equation}\label{eq:illexa2000}
\mathbf{K}(\zeta)=b(\zeta)\,\mathbf{B}+c(\zeta)\,\mathbf{C},
\end{equation}
where
\begin{equation*}
\mathbf{B}=\begin{bmatrix}
2 & -1 & 0\\
-1 & 1 & 0\\
0 & 0 & 0
\end{bmatrix}
\quad\text{and}\quad
\mathbf{C}=\begin{bmatrix}
0 & 0 & 0\\
0 & 1 & -1\\
0 & -1 & 1
\end{bmatrix}.
\end{equation*}

Section \ref{sec:spesol} requires that the random basis is orthogonal with respect to $\mu\sim\mathrm{Uniform}^{\otimes 2}$.
One choice for such a random basis is to consider the tensor product of all Legendre polynomials.
In this example, six random basis functions are utilized to define the discretized version of $\mathscr{Z}$, namely $\mathscr{Z}^{[5]}=\mathrm{span}(\varphi_j)_{j=0}^5$ ($P=5$), where
\begin{gather*}
\varphi_0(\zeta)=1,\quad
\varphi_1(\zeta)=\zeta^1,\quad
\varphi_2(\zeta)=\zeta^2,\\
\varphi_3(\zeta)=\tfrac{3}{2}(\zeta^1)^2-\tfrac{1}{2},\quad
\varphi_4(\zeta)=\zeta^1\zeta^2,\quad
\varphi_5(\zeta)=\tfrac{3}{2}(\zeta^2)^2-\tfrac{1}{2}.
\end{gather*}

By choice, the metric tensor endowed on $\mathbb{R}^3$ is the Euclidean metric tensor, which means that $\eta_{uv}=\delta_{uv}$ (the Kronecker delta).

Substituting Eq.~\eqref{eq:illexa2000} into Eq.~(\ref{eq:spesol700}a) gives the components $K\indices{^u_v^i_j}$ in matrix form:
\begin{equation}\label{eq:illexa3000}
\mathbf{K}\indices{^i_j}=
\frac{\langle\varphi_i,\mathbf{K}\varphi_j\rangle}{\langle\varphi_i,\varphi_i\rangle}=
b\indices{^i_j}\mathbf{B}+c\indices{^i_j}\mathbf{C},
\end{equation}
where $i,j\in\{0,1,\ldots,5\}$, and
\begin{equation}\label{eq:illexa4000}
b\indices{^i_j}=\frac{\langle\varphi_i,b\varphi_j\rangle}{\langle\varphi_i,\varphi_i\rangle}
\quad\text{and}\quad
c\indices{^i_j}=\frac{\langle\varphi_i,c\varphi_j\rangle}{\langle\varphi_i,\varphi_i\rangle}.
\end{equation}

The task now reduces to determining the components $b\indices{^i_j}$, $c\indices{^i_j}$ and $I\indices{^i_{jk}}$ using Eqs.~\eqref{eq:illexa4000} and (\ref{eq:spesol700}b).
This is easy to do and is left to the reader to accomplish.

The entries of $\hat{\mathbf{F}}\in\mathrm{GL}(24,\mathbb{R})$ are determined with Eq.~\eqref{eq:numsol1600} and then mapped using the correspondence set forth in Remark \ref{rmk:numsol2000}.
For example, for $u=1$, $i=2$, $\alpha=2$ and $\beta=3$, we have
\begin{equation*}
\hat{F}\indices{^{ui}_\beta}=
\hat{F}\indices{^{12}_3}=
-\sum_{j=0}^5 I\indices{^2_{3j}}\hat{\phi}^{1j}
\mapsto\text{$(7,4)$-th entry of $\hat{\mathbf{F}}$},
\end{equation*}
\begin{equation*}
\hat{F}\indices{^u_\alpha^i_\beta}=
\hat{F}\indices{^1_2^2_3}=
K\indices{^1_2^2_3}-\delta^1_2\sum_{j=0}^5 I\indices{^2_{3j}}\hat{\lambda}^j
\mapsto\text{$(7,17)$-th entry of $\hat{\mathbf{F}}$},
\end{equation*}
\begin{equation*}
\hat{F}\indices{^i_\beta}=
\hat{F}\indices{^2_3}=
0\mapsto\text{$(21,4)$-th entry of $\hat{\mathbf{F}}$},
\end{equation*}
\begin{equation*}
\hat{F}\indices{_\alpha^i_\beta}=
\hat{F}\indices{_2^2_3}=
-2\sum_{u=1}^3\sum_{j=0}^5 I\indices{^2_{3j}}\delta_{2u}\hat{\phi}^{uj}
\mapsto\text{$(21,17)$-th entry of $\hat{\mathbf{F}}$},
\end{equation*}
where $\delta^1_2=\delta_{21}=\delta_{23}=0$ and $\delta_{22}=1$.

Following a similar approach, the entries of $\hat{\mathbf{H}}\in\mathrm{GL}(24,\mathbb{R})$ are determined with Eq.~\eqref{eq:numsol2200} and then mapped using the correspondence set forth in Remark \ref{rmk:numsol3000}.
Alternatively, $\hat{\mathbf{H}}$ can be computed with Eq.~\eqref{eq:numsol1600} provided that $\hat{x}$ is updated with $\bar{x}$ (see Remark \ref{rmk:numsol6000}).
The entries of $\hat{\mathbf{X}}\in\mathrm{GL}(24,\mathbb{R})$ are determined with Eq.~\eqref{eq:numsol1650}.
The entries of $\bar{\mathbf{x}},\hat{\mathbf{x}},\hat{\mathbf{f}}\in\mathbb{R}^{24}$ are obtained with Eqs.~\eqref{eq:numsol2300} and \eqref{eq:numsol2400}.

The algorithm outlined below provides the numerical solution to Eq.~\eqref{eq:illexa1000}.
\begin{enumerate}
\item \emph{Initial guess:} Choose an initial guess $\hat{\mathbf{x}}$ close to the expected solution.
This can be achieved by conducting a Monte Carlo simulation with a few number of realizations to position the initial guess sufficiently close to the solution.
\item \emph{Iteration:} Compute $\hat{\mathbf{F}}$, $\hat{\mathbf{X}}$, $\hat{\mathbf{H}}$, $\hat{\mathbf{L}}$ and $\hat{\mathbf{f}}$. 
Then, update the guess using:
\begin{equation*}
\mathbf{x}=\hat{\mathbf{x}}-\hat{\mathbf{L}}\hat{\mathbf{f}}.
\end{equation*}
This follows from Eq.~\eqref{eq:numsol1150} with $\hat{\mathbf{L}}=\hat{\mathbf{H}}^{-1}$.
(We remark that for Newton's method, $\hat{\mathbf{L}}$ is taken as $\hat{\mathbf{X}}$.)
\item \emph{Convergence check:} Repeat step 2 until the change in magnitude of $\mathbf{x}$ is smaller than a pre-specified tolerance $\epsilon$ (that is, $\lVert\hat{\mathbf{x}}-\mathbf{x}\rVert_2\equiv\lVert\hat{\mathbf{L}}\hat{\mathbf{f}}\rVert_2<\epsilon$), or until the magnitude of $\hat{\mathbf{f}}$ is smaller than $\epsilon$ (that is, $\lVert\hat{\mathbf{f}}\rVert_2\equiv\lVert\hat{\mathbf{H}}(\hat{\mathbf{x}}-\mathbf{x})\rVert_2<\epsilon$).
In this work, we use the latter check, where $\lVert\,\cdot\,\rVert_2$ denotes the Euclidean norm on $\mathbb{R}^Q$ ($Q=24$ for $R=3$ and $P=5$).
\end{enumerate}

Implementing the algorithm with $\epsilon=10^{-12}$ and $P=5$ yields
\begin{equation*}
\mathbf{E}[\lambda]=
\begin{cases}
3.73451918626129 & \text{for mode 1}\\
25.4620874810967 & \text{for mode 2}\\
61.0547534054944 & \text{for mode 3}
\end{cases}
\qquad\text{and}\qquad
\mathrm{Var}[\lambda]=
\begin{cases}
0.28824486249498 & \text{for mode 1}\\
38.3590753966161 & \text{for mode 2}\\
70.6920253818825 & \text{for mode 3},
\end{cases}
\end{equation*}
which follow from using Eqs.~\eqref{eq:promom300} and \eqref{eq:promom400}, respectively.

\begin{figure}
\centering
\includegraphics[scale=0.8]{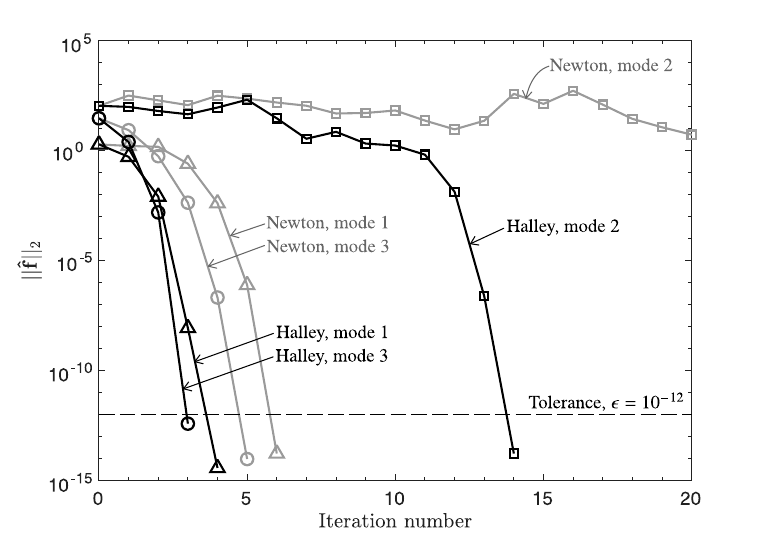}
\caption{Convergence of Halley's method and Newton's method}
\label{fig:illexa1000}
\end{figure}

Fig.~\ref{fig:illexa1000} compares the convergence of Newton's method and Halley's
method.
It is observed that Halley's method not only exhibits a faster convergence rate than Newton's method, but also converges in cases where Newton's method may not, highlighting the robustness and efficiency of Halley's approach.

\section{A case study}\label{sec:casstu}

\begin{figure}
\centering
\includegraphics[scale=1.0]{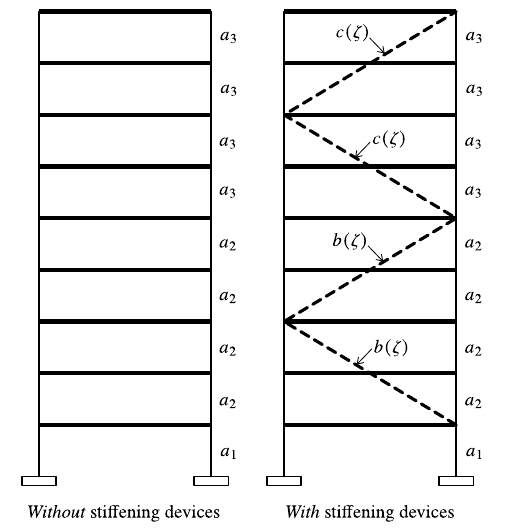}
\caption{Nine-story building under study}
\label{fig:casstu1000}
\end{figure}

To demonstrate our method with a numerical example, we consider a 9-story office building located in a hurricane-prone region. 
The building (Fig.~\ref{fig:casstu1000}) is designed using a combination of frame elements and stiffening devices to enhance its dynamic performance.
The primary goal of this design is to assess whether the incorporation of stiffening devices along the height of the building can effectively improve the structural integrity of the building, increase occupant comfort, and distribute wind loads more evenly throughout the structure.

A linearized version of the differential equation governing the motion of the 9-story building, without including damping effects or applied forces, is:
\begin{equation}
\mathbf{M}\ddot{\mathbf{u}}+\mathbf{R}\mathbf{u}=\mathbf{0},
\end{equation}
where $\mathbf{R}$ is the stiffness matrix, $\mathbf{M}$ is the mass matrix, $\mathbf{u}$ is the displacement vector, $\ddot{\mathbf{u}}$ is the acceleration vector, and $\mathbf{0}$ is the zero vector.

For the purposes of this example, we consider nine translational (horizontal) degrees of freedom, one for each story, and take the mass matrix as $\mathbf{M}=m\mathbf{I}$, where $m$ is the inter-story mass of the building in units of Mg, and $\mathbf{I}$ is the 9-by-9 identity matrix.
As for the stiffness matrix, we take it as
\begin{equation}
\mathbf{R}(\zeta)=\mathbf{A}+b(\zeta)\,\mathbf{B}+c(\zeta)\,\mathbf{C},
\end{equation}
where $b,c$ are real-valued, uniformly-distributed random variables in units of kN/m, representing the lateral stiffnesses provided by the lower-level and upper-level stiffening devices, respectively, and $\mathbf{R},\mathbf{A},\mathbf{B},\mathbf{C}$ are 9-by-9 real symmetric matrices with the entries of $\mathbf{R}$ and $\mathbf{A}$ in units of kN/m and the entries of $\mathbf{B}$ and $\mathbf{C}$ unitless.
Specifically, $b(\zeta)=(50+7.5\zeta^1)\times10^3$, $c(\zeta)=(25+3.75\zeta^2)\times10^3$, $\zeta=(\zeta^1,\zeta^2)\sim\mathrm{Uniform}^{\otimes 2}$, and the random domain is $\mathfrak{Z}=[-1,1]^2$.
The nonzero lower-triangular entries of $\mathbf{A}$ are: 
\begin{gather*}
A\indices{^1_1}=a_1+a_2,\quad 
A\indices{^2_2}=A\indices{^3_3}=A\indices{^4_4}=2a_2,\quad
A\indices{^5_5}=a_2+a_3,\quad 
A\indices{^6_6}=A\indices{^7_7}=A\indices{^8_8}=2a_3,\quad 
A\indices{^9_9}=a_3,\\
A\indices{^2_1}=A\indices{^3_2}=A\indices{^4_3}=A\indices{^5_4}=-a_2,\quad 
A\indices{^6_5}=A\indices{^7_6}=A\indices{^8_7}=A\indices{^9_8}=-a_3,
\end{gather*} 
where $a_1=170\times10^3$, $a_2=120\times10^3$ and $a_3=75\times10^3$ represent the lateral stiffnesses of the first story, the second-to-fifth story, and the sixth-to-ninth story, respectively.
These stiffnesses are expressed in units of kN/m.
The nonzero lower-triangular entries of $\mathbf{B}$ and $\mathbf{C}$ are: 
\begin{gather*}
B\indices{^1_1}=B\indices{^5_5}=C\indices{^5_5}=C\indices{^9_9}=1,\quad 
B\indices{^3_3}=C\indices{^7_7}=2,\\ 
B\indices{^3_1}=B\indices{^5_3}=C\indices{^7_5}=C\indices{^9_7}=-1.
\end{gather*}

If the mean values of $b$ and $c$ are used and $m=100$ Mg, the fundamental period of the building without the stiffening devices is 1.12 s, whereas with the stiffening devices it is 0.88 s, indicating a significant reduction of more than 20\% in the fundamental period of the building when the stiffening devices are employed.

We now solve the stochastic eigenvalue problem stated in Eq.~\eqref{eq:prosta100} by setting  $\mathbf{K}=\mathbf{R}$.
Thus, the eigenvalues are expressed in units of $\mathrm{Mg}/\mathrm{s}^2$, and the period of vibration of the building is computed with the expression: $T=2\pi\sqrt{m/\lambda}$.
The solution is numerically obtained using $P=5$ and $\epsilon=10^{-8}$.
For reference, the nonrandom eigenvalues are listed in the second column of Table \ref{tab:casstu1000} using the mean values of $b$ and $c$; that is, $b=50\times10^3$ and $c=25\times10^3$ (both in kN/m).

\begin{figure}
\centering
\includegraphics[scale=0.8]{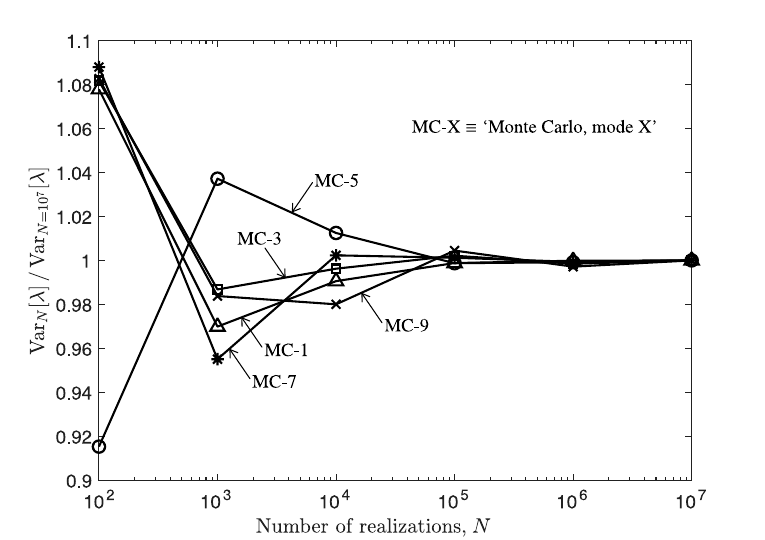}
\caption{Variance convergence of Monte Carlo simulations for the 1st, 3rd, 5th, 7th and 9th modes (variances are normalized relative to the variance obtained from $N=10^7$ realizations)}
\label{fig:casstu2000}
\end{figure}

The solution outline is as follows.
First, six Monte Carlo simulations are employed to estimate the chaos coefficients.
The coefficients are computed using different sample sizes to establish the number of realizations needed to get a reliable estimate.
Fig.~\ref{fig:casstu2000} illustrates the variance convergence of Monte Carlo simulations as a function of the sample size for the system's odd modes.
When $N$ is relatively small, significant deviations from the true variance (Monte Carlo simulation with $N=10^7$ realizations) are observed, as expected.
However, as $N$ increases, the variance estimates progressively converge to the true value thanks to the law of large numbers.
This is clearly illustrated in the figure, where the computed variances become more tightly concentrated around the true value for large $N$.

\begin{figure}
\centering
\includegraphics[scale=0.7]{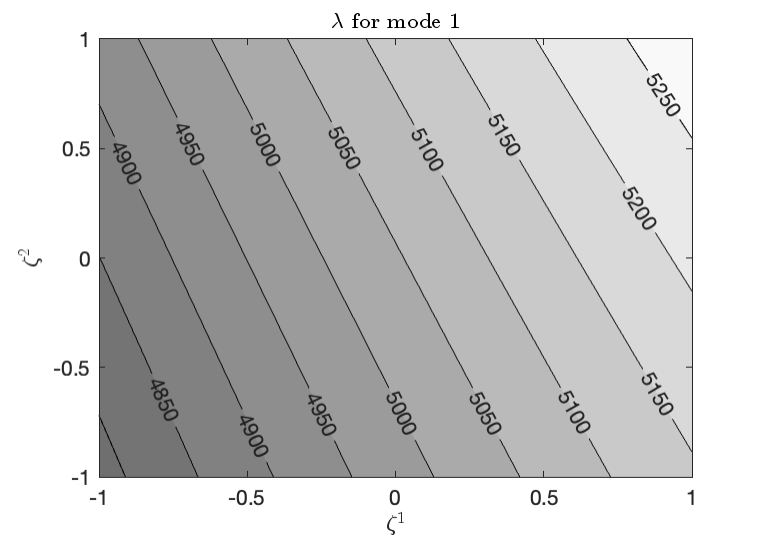}\\
\includegraphics[scale=0.7]{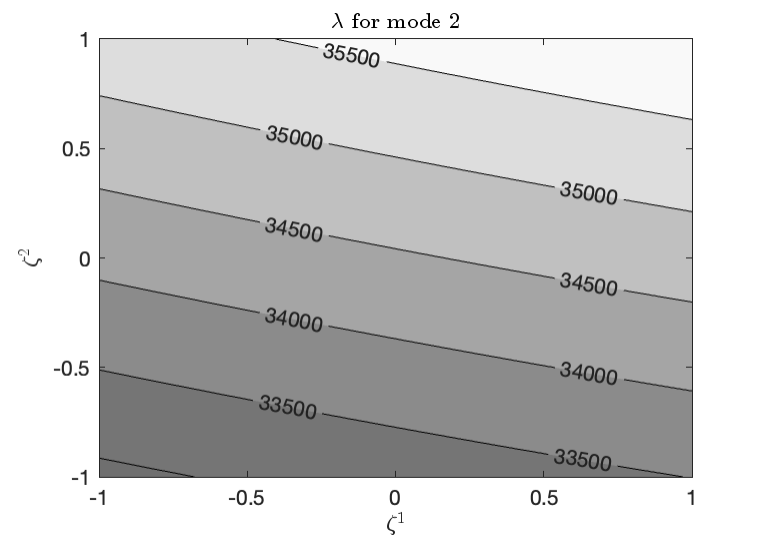}\\
\includegraphics[scale=0.7]{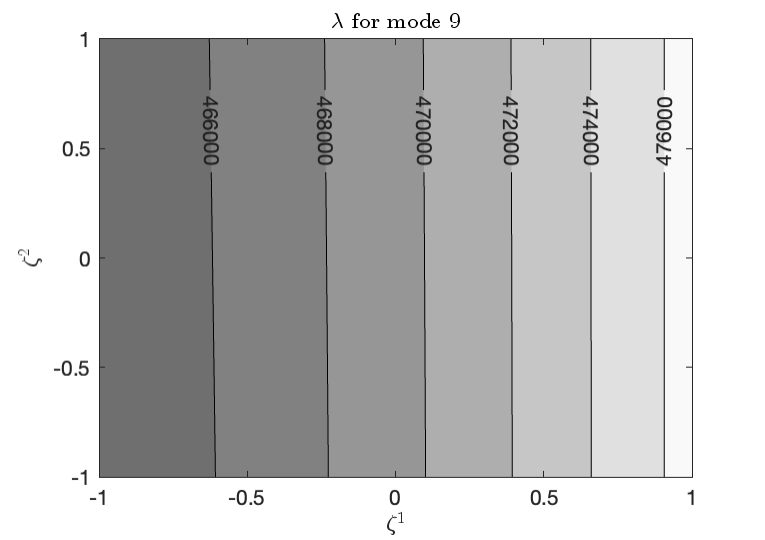}
\caption{Plots of the 1st, 2nd and 9th eigenvalues as functions of $\zeta=(\zeta^1,\zeta^2)$ using Halley's method}
\label{fig:casstu3000}
\end{figure}

Fig.~\ref{fig:casstu3000} shows the relationship between the first, second and ninth eigenvalues and the random domain using Halley's method.
The random domain is represented by a square region in the $\zeta^1\zeta^2$-plane, ranging from $-1$ to $1$ in both directions, while the eigenvalues are represented by contour lines.
This visualization is essential for understanding how the eigenvalues respond to variations in the random inputs.
The results show that an increase in the values of $\zeta^1$ and $\zeta^2$ leads to a gradual increase in $\lambda$, specifically for the first-to-fourth and seventh modes.
In contrast, for the eighth and ninth modes, $\lambda$ increases primarily as a function of $\zeta^1$, while for the fifth and sixth modes, it increases mainly as a function of $\zeta^2$.
These findings provide valuable information on the stability and sensitivity of the system over the random domain.

\begin{figure}
\centering
\includegraphics[scale=0.8]{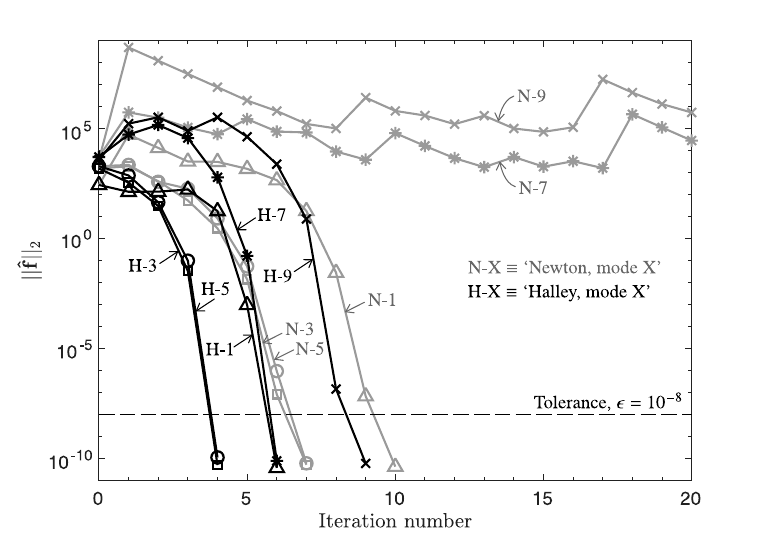}
\caption{Convergence of Halley's method and Newton's method for the 1st, 3rd, 5th, 7th and 9th modes}
\label{fig:casstu4000}
\end{figure}

Fig.~\ref{fig:casstu4000} illustrates the comparative performance of Newton's and Halley's methods as a function of the iteration number.
In this context, iteration number zero represents the initial errors before the iterative scheme begins.
The dashed line indicates the pre-specified tolerance of $\epsilon=10^{-8}$.
The results demonstrate that Halley's method consistently converges faster than Newton's method.
Specifically, Halley's method achieves the tolerance level in fewer iterations, with
some modes converging in as few as 4 iterations.
This efficiency suggests that the higher-order correction terms of Halley's method offer a significant advantage in rapidly reducing errors.
In contrast, Newton's method, while still converging effectively for most of the modes, requires more iterations to achieve the same level of tolerance.

Finally, Table \ref{tab:casstu1000} presents the results of a comparative performance analysis based on a Monte Carlo simulation with $N=10^6$ realizations, and Newton's and Halley's method with 100 iterations allowed.
The table shows that, while all three methods provide comparable eigenvalue estimates, Halley's method consistently achieved the more reliable estimates among the three methods.
In particular, it is noted that for the seventh mode, Newton's method fails to converge within the 100 iterations allotted, while for the ninth mode, important discrepancies in the mean and variance values are obtained after Newton's method identified an alternative, though still incorrect, solution.

\begin{landscape}
\begin{table}
\caption{Probability moments of all eigenvalues as computed by a Monte Carlo simulation with $N=10^6$ realizations, Halley's method and Newton's method}
\label{tab:casstu1000}
\small
\begin{tabular}{cccccccccc}
\toprule
Mode & $\lambda\in\mathbb{R}^+$ & \multicolumn{3}{c}{Mean: $\mathbf{E}[\lambda]\in\mathbb{R}^+$} & \multicolumn{3}{c}{Standard deviation: $(\mathrm{Var}[\lambda])^{1/2}\in\mathbb{R}^+_0$} & \multicolumn{2}{c}{Number of iterations} \\
\cmidrule(lr){3-5} \cmidrule(lr){6-8} \cmidrule(lr){9-10}
   & Nonrandom  & Monte Carlo & Halley & Newton & Monte Carlo & Halley & Newton & Halley & Newton \\
\midrule
1   & 5041.69   & 5038.48   & 5040.82   & 5040.82         & 112.14  & 112.06  & 112.06          & 6 & 10         \\
2   & 34442.52  & 34434.34  & 34435.22  & 34435.22        & 714.66  & 719.74  & 719.74          & 8 & 9          \\
3   & 93274.43  & 93214.44  & 93243.65  & 93243.65        & 1331.34 & 1336.97 & 1336.97         & 4 & 7          \\
4   & 160320.79 & 160249.62 & 160291.86 & 160291.86       & 1983.91 & 1969.54 & 1969.54         & 4 & 7          \\
5   & 192832.06 & 192864.51 & 192851.29 & 192851.29       & 1275.12 & 1271.53 & 1271.53         & 4 & 7          \\
6   & 284340.38 & 284387.06 & 284358.30 & 284358.30       & 2500.77 & 2498.28 & 2498.28         & 4 & 7          \\
7   & 376567.59 & 376166.88 & 376286.25 & {\it 373622.34} & 6973.61 & 7155.76 & {\it 43598.02}  & 6 & {\it 100} $^\dag$ \\
8   & 413799.29 & 413866.90 & 413864.51 & 413864.51       & 4634.62 & 4643.04 & 4643.04         & 4 & 7          \\
9   & 469381.25 & 469763.46 & 469779.80 & {\it 403882.22} & 3605.96 & 3606.22 & {\it 148806.99} & 9 & {\it 33} $^\ddag$  \\
\bottomrule\\[-2ex]
\multicolumn{10}{l}{$^\dag$ The solution did not converge after 100 iterations.}
\\
\multicolumn{10}{l}{$^\ddag$ The solution converged to a different root after 33 iterations.}
\end{tabular}
\end{table}
\end{landscape}

\section{Conclusion}\label{sec:conclu}

We have presented a novel numerical scheme to solve the stochastic eigenproblem via the spectral-chaos approach, using Halley's method as the root-finding algorithm for the resulting set of deterministic quadratic equations displayed in Eq.~\eqref{eq:spesol600}.
Halley's method was chosen not only for its cubic convergence properties but also because a higher-order Householder method would not provide additional improvements in the convergence rate, given the quadratic nature of the resulting system of equations.
Our method clearly marks a substantial improvement over the approach of Ghanem and Ghosh \cite{ghanem2007efficient}, because it achieves a higher polynomial order of convergence and has the ability to accelerate the convergence rate when the eigenvector components are nearly known.
Moreover, to make the solution tractable, the stochastic eigenvalue problem was formulated in tensor notation to significantly reduce the complexity of the solution.
To the best of the authors' knowledge, this approach had not been explored before as a compelling alternative to solving the stochastic eigenvalue problem.

Both the illustrative example and the case study demonstrated cases where Halley's method converged while Newton's method did not.
The comparative performance analyses revealed that Halley's method converges faster and more reliably than Newton's method.
Additionally, Halley's method consistently delivered the most accurate results among the three methods.
Hence, Halley's method is recommended for applications that require higher reliability and accuracy.

The algorithmic difference between Newton's method and Halley's method lies in the need to compute two extra objects for Halley's case: $\hat{\mathbf{H}}$ and $\hat{\mathbf{L}}$.
This additional computation increases the computational cost of Halley's method compared to Newton's method, resulting in a slower execution in each iteration.
However, the potential slowdown of Halley's method can be offset by its fast convergence rate, especially in cases where Newton's method struggles to find the solution due to an inadequate initial guess.
In Section \ref{sec:errana}, it was shown that Halley's method has the ability to accelerate the solution process when the eigenvector components are nearly known, an advantage not shared by Newton's method.
Section \ref{sec:comcos} demonstrated that both Halley's and Newton's methods exhibit the same cubic asymptotic time complexity, indicating thereby that Halley's method scales just as efficiently as Newton's method in both high-dimensional spatial and random function spaces.

The authors believe, nonetheless, that further efforts are necessary to develop a numerical scheme that effectively aids in solving the stochastic eigenvalue problem with greater efficiency. 
This can be accomplished by integrating advanced techniques from the current literature into Newton-based methods or by exploring innovative new approaches.

\section*{Author contributions}
H.E.~developed the conceptual framework and methodology of the proposed method (Sections 2 to 8).
H.E.~and K.O.I.~presented and illustrated the proposed method, and carried out the implementation and validation of the discussed methods.
G.L.~provided the funding resources for K.O.I., and supervised the findings of this work.
All authors discussed the results and contributed to the final manuscript.

\section*{Acknowledgements}
This work was supported by the National Science Foundation (DMS-2053746, DMS-2134209, ECCS-2328241, CBET-2347401 and OAC-2311848), and U.S.~Department of Energy (DOE) Office of Science Advanced Scientific Computing Research Program DE-SC0023161, and DOE-Fusion Energy Science under grant number: DE-SC0024583.
The first author wishes to express sincere gratitude to Estructuras del Norte, a Colombian structural engineering firm, for their financial support.
Part of this work was developed while the first author was an Assistant Professor in the Department of Civil and Environmental Engineering at Universidad de la Costa---their support is gratefully acknowledged.

\section*{Financial disclosure}
There are no financial conflicts of interest to disclose.

\section*{Conflict of interest}
The authors declare no potential conflict of interests.

\bibliographystyle{amsplain}

\bibliography{References}

\newpage

\section*{Appendix A}

In this section, we derive the derivative of $F$ with respect to $f=(g,h)$.
Specifically, the eight components produced by
\begin{equation}
\frac{\partial}{\partial f^s}\left(\frac{\partial f^n}{\partial x^m}\right)\bigg|_{f(\hat{x})}\equiv
\frac{\partial F\indices{^n_m}}{\partial f^s}\bigg|_{f(\hat{x})}\tag{\ref{eq:numsol1700}}
\end{equation}
are presented.
The indices satisfy $n\in(ui,i)$, $m\in(\beta,\alpha\beta)$ and $s\in(\gamma\rho,\rho)$, where $u,\alpha,\gamma\in\{1,2,\ldots,R\}$ and $i,j,\beta,\rho\in\{0,1,\ldots,P\}$.

Thus, differentiating Eq.~\eqref{eq:numsol1600} with respect to $f^s=(g^{\gamma\rho},h^\rho)$ yields:
\begin{subequations}
\begin{equation}
\frac{\partial F\indices{^{ui}_\beta}}{\partial g^{\gamma\rho}}\bigg|_{f(\hat{x})}=
-I\indices{^i_{\beta j}}\frac{\partial\phi^{uj}}{\partial g^{\gamma\rho}}\bigg|_{f(\hat{x})}=
-I\indices{^i_{\beta j}}\hat{X}\indices{^u_\gamma^j_\rho}.
\end{equation}
\begin{equation}
\frac{\partial F\indices{^{ui}_\beta}}{\partial h^\rho}\bigg|_{f(\hat{x})}=
-I\indices{^i_{\beta j}}\frac{\partial\phi^{uj}}{\partial h^\rho}\bigg|_{f(\hat{x})}=
-I\indices{^i_{\beta j}}\hat{X}\indices{^{uj}_\rho}.
\end{equation}
\begin{equation}
\frac{\partial F\indices{^u_\alpha^i_\beta}}{\partial g^{\gamma\rho}}\bigg|_{f(\hat{x})}=
-I\indices{^i_{\beta j}}\delta^u_\alpha\frac{\partial\lambda^j}{\partial g^{\gamma\rho}}\bigg|_{f(\hat{x})}=
-I\indices{^i_{\beta j}}\delta^u_\alpha\hat{X}\indices{_\gamma^j_\rho}.
\end{equation}
\begin{equation}
\frac{\partial F\indices{^u_\alpha^i_\beta}}{\partial h^\rho}\bigg|_{f(\hat{x})}=
-I\indices{^i_{\beta j}}\delta^u_\alpha\frac{\partial\lambda^j}{\partial h^\rho}\bigg|_{f(\hat{x})}=
-I\indices{^i_{\beta j}}\delta^u_\alpha\hat{X}\indices{^j_\rho}.
\end{equation}
\begin{equation}
\frac{\partial F\indices{^i_\beta}}{\partial g^{\gamma\rho}}\bigg|_{f(\hat{x})}=
\frac{\partial F\indices{^i_\beta}}{\partial h^\rho}\bigg|_{f(\hat{x})}=0.
\end{equation}
\begin{equation}
\frac{\partial F\indices{_\alpha^i_\beta}}{\partial g^{\gamma\rho}}\bigg|_{f(\hat{x})}=
2I\indices{^i_{\beta j}}\eta_{\alpha u}\frac{\partial\phi^{uj}}{\partial g^{\gamma\rho}}\bigg|_{f(\hat{x})}=
2I\indices{^i_{\beta j}}\eta_{\alpha u}\hat{X}\indices{^u_\gamma^j_\rho}.
\end{equation}
\begin{equation}
\frac{\partial F\indices{_\alpha^i_\beta}}{\partial h^\rho}\bigg|_{f(\hat{x})}=
2I\indices{^i_{\beta j}}\eta_{\alpha u}\frac{\partial\phi^{uj}}{\partial h^\rho}\bigg|_{f(\hat{x})}=
2I\indices{^i_{\beta j}}\eta_{\alpha u}\hat{X}\indices{^{uj}_\rho}.
\end{equation}
\end{subequations}

\section*{Appendix B}

In this section, we derive the derivative of $F$ with respect to $x=(\lambda,\phi)$:
\begin{equation}
\frac{\partial^2 f^s}{\partial x^a\partial x^b}\bigg|_y\equiv\frac{\partial F\indices{^s_b}}{\partial x^a}\bigg|_y.\tag{\ref{eq:errana13000}}
\end{equation}

Similarly to Appendix A, this expression yields eight components.
The indices satisfy $s\in(ui,i)$, $b\in(\beta,\alpha\beta)$ and $a\in(\rho,\gamma\rho)$, where $u,\alpha,\gamma\in\{1,2,\ldots,R\}$ and $i,j,\beta,\rho\in\{0,1,\ldots,P\}$.
Hence, differentiating Eq.~\eqref{eq:numsol1600} with respect to $x^s=(\lambda^\rho,\phi^{\gamma\rho})$ produces:
\begin{subequations}
\begin{equation}
\frac{\partial F\indices{^{ui}_\beta}}{\partial\lambda^\rho}\bigg|_y
=\frac{\partial F\indices{^u_\alpha^i_\beta}}{\partial\phi^{\gamma\rho}}\bigg|_y
=\frac{\partial F\indices{^i_\beta}}{\partial\lambda^\rho}\bigg|_y
=\frac{\partial F\indices{^i_\beta}}{\partial\phi^{\gamma\rho}}\bigg|_y
=\frac{\partial F\indices{_\alpha^i_\beta}}{\partial\lambda^\rho}\bigg|_y
=0.
\end{equation}
\begin{equation}
\frac{\partial F\indices{^{ui}_\beta}}{\partial\phi^{\gamma\rho}}\bigg|_y=
-I\indices{^i_{\beta j}}\frac{\partial\phi^{uj}}{\partial\phi^{\gamma\rho}}\bigg|_y=
-I\indices{^i_{\beta j}}\delta^u_\gamma\delta^j_\rho=
-I\indices{^i_{\beta\rho}}\delta^u_\gamma.
\end{equation}
\begin{equation}
\frac{\partial F\indices{^u_\alpha^i_\beta}}{\partial\lambda^\rho}\bigg|_y=
-I\indices{^i_{\beta j}}\frac{\partial\lambda^j}{\partial\lambda^\rho}\bigg|_y\delta^u_\alpha=
-I\indices{^i_{\beta j}}\delta^j_\rho\delta^u_\alpha=
-I\indices{^i_{\beta\rho}}\delta^u_\alpha.
\end{equation}
\begin{equation}
\frac{\partial F\indices{_\alpha^i_\beta}}{\partial\phi^{\gamma\rho}}\bigg|_y=
2I\indices{^i_{\beta j}}\eta_{\alpha u}\frac{\partial\phi^{uj}}{\partial\phi^{\gamma\rho}}\bigg|_y=
2I\indices{^i_{\beta j}}\eta_{\alpha u}\delta^u_\gamma\delta^j_\rho=
2I\indices{^i_{\beta\rho}}\eta_{\alpha\gamma}.
\end{equation}
\end{subequations}

\section*{Appendix C}

In this section, we derive the individual components of the composite object:
\begin{equation}
\check{D}\indices{^n_{abc}}:=\frac{\partial^2 f^s}{\partial x^a\partial x^b}\bigg|_y\frac{\partial}{\partial f^s}\bigg(\frac{\partial f^n}{\partial x^c}\bigg)\bigg|_{f(y)}\tag{\ref{eq:errana15000}}.
\end{equation}
The indices satisfy $n\in(vj,j)$, $c\in(\sigma,\omega\sigma)$, $b\in(\beta,\alpha\beta)$, $a\in(\rho,\gamma\rho)$ and $s\in(ui,i)$, where $u,v,w,\alpha,\gamma,\omega\in\{1,2,\ldots,R\}$ and $i,j,k,\beta,\rho,\sigma\in\{0,1,\ldots,P\}$.

\begin{subequations}\label{eq:errana16000}
From Appendices A and B, it should be clear that:
\begin{equation}
\check{D}\indices{^{(vj)}_{\rho\beta\sigma}}
=\check{D}\indices{^{(vj)}_{\rho\beta(\omega\sigma)}}
=\check{D}\indices{^j_{\rho\beta\sigma}}
=\check{D}\indices{^j_{(\gamma\rho)\beta\sigma}}
=\check{D}\indices{^j_{\rho(\alpha\beta)\sigma}}
=\check{D}\indices{^j_{(\gamma\rho)(\alpha\beta)\sigma}}
=\check{D}\indices{^j_{\rho\beta(\omega\sigma)}}
=0.
\end{equation}
The other components of $\check{D}\indices{^n_{abc}}$ are as follows:
\begin{align}
\check{D}\indices{^{(vj)}_{(\gamma\rho)\beta\sigma}}&=
\frac{\partial F\indices{^{ui}_\beta}}{\partial\phi^{\gamma\rho}}\bigg|_y\frac{\partial F\indices{^{vj}_\sigma}}{\partial g^{ui}}\bigg|_{f(y)}+\frac{\partial F\indices{^i_\beta}}{\partial\phi^{\gamma\rho}}\bigg|_y\frac{\partial F\indices{^{vj}_\sigma}}{\partial h^i}\bigg|_{f(y)}\notag\\
&=(-I\indices{^i_{\beta\rho}}\delta^u_\gamma)(-I\indices{^j_{\sigma k}}\check{X}\indices{^v_u^k_i})+0
=I\indices{^i_{\beta\rho}}I\indices{^j_{\sigma k}}\check{X}\indices{^v_\gamma^k_i}.
\end{align}
\begin{align}
\check{D}\indices{^{(vj)}_{\rho(\alpha\beta)\sigma}}&=
\frac{\partial F\indices{^u_\alpha^i_\beta}}{\partial\lambda^\rho}\bigg|_y\frac{\partial F\indices{^{vj}_\sigma}}{\partial g^{ui}}\bigg|_{f(y)}+\frac{\partial F\indices{_\alpha^i_\beta}}{\partial\lambda^\rho}\bigg|_y\frac{\partial F\indices{^{vj}_\sigma}}{\partial h^i}\bigg|_{f(y)}\notag\\
&=(-I\indices{^i_{\beta\rho}}\delta^u_\alpha)(-I\indices{^j_{\sigma k}}\check{X}\indices{^v_u^k_i})+0
=I\indices{^i_{\beta\rho}}I\indices{^j_{\sigma k}}\check{X}\indices{^v_\alpha^k_i}.
\end{align}
\begin{align}
\check{D}\indices{^{(vj)}_{(\gamma\rho)(\alpha\beta)\sigma}}&=
\frac{\partial F\indices{^u_\alpha^i_\beta}}{\partial\phi^{\gamma\rho}}\bigg|_y\frac{\partial F\indices{^{vj}_\sigma}}{\partial g^{ui}}\bigg|_{f(y)}+\frac{\partial F\indices{_\alpha^i_\beta}}{\partial\phi^{\gamma\rho}}\bigg|_y\frac{\partial F\indices{^{vj}_\sigma}}{\partial h^i}\bigg|_{f(y)}\notag\\
&=0+(2I\indices{^i_{\beta\rho}}\eta_{\alpha\gamma})(-I\indices{^j_{\sigma k}}\check{X}\indices{^{vk}_i})
=-2I\indices{^i_{\beta\rho}}I\indices{^j_{\sigma k}}\eta_{\alpha\gamma}\check{X}\indices{^{vk}_i}.
\end{align}
\begin{align}
\check{D}\indices{^{(vj)}_{(\gamma\rho)\beta(\omega\sigma)}}&=
\frac{\partial F\indices{^{ui}_\beta}}{\partial\phi^{\gamma\rho}}\bigg|_y\frac{\partial F\indices{^v_\omega^j_\sigma}}{\partial g^{ui}}\bigg|_{f(y)}+\frac{\partial F\indices{^i_\beta}}{\partial\phi^{\gamma\rho}}\bigg|_y\frac{\partial F\indices{^v_\omega^j_\sigma}}{\partial h^i}\bigg|_{f(y)}\notag\\
&=(-I\indices{^i_{\beta\rho}}\delta^u_\gamma)(-I\indices{^j_{\sigma k}}\delta^v_\omega\check{X}\indices{_u^k_i})+0
=I\indices{^i_{\beta\rho}}I\indices{^j_{\sigma k}}\delta^v_\omega\check{X}\indices{_\gamma^k_i}.
\end{align}
\begin{align}
\check{D}\indices{^{(vj)}_{\rho(\alpha\beta)(\omega\sigma)}}&=
\frac{\partial F\indices{^u_\alpha^i_\beta}}{\partial\lambda^\rho}\bigg|_y\frac{\partial F\indices{^v_\omega^j_\sigma}}{\partial g^{ui}}\bigg|_{f(y)}+\frac{\partial F\indices{_\alpha^i_\beta}}{\partial\lambda^\rho}\bigg|_y\frac{\partial F\indices{^v_\omega^j_\sigma}}{\partial h^i}\bigg|_{f(y)}\notag\\
&=(-I\indices{^i_{\beta\rho}}\delta^u_\alpha)(-I\indices{^j_{\sigma k}}\delta^v_\omega\check{X}\indices{_u^k_i})+0
=I\indices{^i_{\beta\rho}}I\indices{^j_{\sigma k}}\delta^v_\omega\check{X}\indices{_\alpha^k_i}.
\end{align}
\begin{align}
\check{D}\indices{^{(vj)}_{(\gamma\rho)(\alpha\beta)(\omega\sigma)}}&=
\frac{\partial F\indices{^u_\alpha^i_\beta}}{\partial\phi^{\gamma\rho}}\bigg|_y\frac{\partial F\indices{^v_\omega^j_\sigma}}{\partial g^{ui}}\bigg|_{f(y)}+\frac{\partial F\indices{_\alpha^i_\beta}}{\partial\phi^{\gamma\rho}}\bigg|_y\frac{\partial F\indices{^v_\omega^j_\sigma}}{\partial h^i}\bigg|_{f(y)}\notag\\
&=0+(2I\indices{^i_{\beta\rho}}\eta_{\alpha\gamma})(-I\indices{^j_{\sigma k}}\delta^v_\omega\check{X}\indices{^k_i})
=-2I\indices{^i_{\beta\rho}}I\indices{^j_{\sigma k}}\eta_{\alpha\gamma}\delta^v_\omega\check{X}\indices{^k_i}.
\end{align}
\begin{align}
\check{D}\indices{^j_{(\gamma\rho)\beta(\omega\sigma)}}&=
\frac{\partial F\indices{^{ui}_\beta}}{\partial\phi^{\gamma\rho}}\bigg|_y\frac{\partial F\indices{_\omega^j_\sigma}}{\partial g^{ui}}\bigg|_{f(y)}+\frac{\partial F\indices{^i_\beta}}{\partial\phi^{\gamma\rho}}\bigg|_y\frac{\partial F\indices{_\omega^j_\sigma}}{\partial h^i}\bigg|_{f(y)}\notag\\
&=(-I\indices{^i_{\beta\rho}}\delta^u_\gamma)(2I\indices{^j_{\sigma k}}\eta_{\omega w}\check{X}\indices{^w_u^k_i})+0
=-2I\indices{^i_{\beta\rho}}I\indices{^j_{\sigma k}}\eta_{\omega w}\check{X}\indices{^w_\gamma^k_i}.
\end{align}
\begin{align}
\check{D}\indices{^j_{\rho(\alpha\beta)(\omega\sigma)}}&=
\frac{\partial F\indices{^u_\alpha^i_\beta}}{\partial\lambda^\rho}\bigg|_y\frac{\partial F\indices{_\omega^j_\sigma}}{\partial g^{ui}}\bigg|_{f(y)}+\frac{\partial F\indices{_\alpha^i_\beta}}{\partial\lambda^\rho}\bigg|_y\frac{\partial F\indices{_\omega^j_\sigma}}{\partial h^i}\bigg|_{f(y)}\notag\\
&=(-I\indices{^i_{\beta\rho}}\delta^u_\alpha)(2I\indices{^j_{\sigma k}}\eta_{\omega w}\check{X}\indices{^w_u^k_i})+0
=-2I\indices{^i_{\beta\rho}}I\indices{^j_{\sigma k}}\eta_{\omega w}\check{X}\indices{^w_\alpha^k_i}.
\end{align}
\begin{align}
\check{D}\indices{^j_{(\gamma\rho)(\alpha\beta)(\omega\sigma)}}&=
\frac{\partial F\indices{^u_\alpha^i_\beta}}{\partial\phi^{\gamma\rho}}\bigg|_y\frac{\partial F\indices{_\omega^j_\sigma}}{\partial g^{ui}}\bigg|_{f(y)}+\frac{\partial F\indices{_\alpha^i_\beta}}{\partial\phi^{\gamma\rho}}\bigg|_y\frac{\partial F\indices{_\omega^j_\sigma}}{\partial h^i}\bigg|_{f(y)}\notag\\
&=0+(2I\indices{^i_{\beta\rho}}\eta_{\alpha\gamma})(2I\indices{^j_{\sigma k}}\eta_{\omega w}\check{X}\indices{^{wk}_i})
=4I\indices{^i_{\beta\rho}}I\indices{^j_{\sigma k}}\eta_{\alpha\gamma}\eta_{\omega w}\check{X}\indices{^{wk}_i}.
\end{align}
\end{subequations}

\end{document}